\documentclass[a4paper,oneside]{amsart}
\usepackage[english]{babel}
\usepackage{amsmath, amssymb, amsthm, amscd}
\usepackage{enumerate}
\usepackage{palatino}
\usepackage{mathpazo}
\usepackage{paralist}
\usepackage[a4paper]{geometry}
\usepackage{url}
\usepackage{hyperref}
\usepackage{tikz}
\usepackage{tikz-cd}
\usepackage{comment}
\usepackage{enumitem}

\theoremstyle{plain}
\newtheorem{theorem}{Theorem}[section]

\newtheorem{lemma}[theorem]{Lemma}
\newtheorem{cor}[theorem]{Corollary}  

\newtheorem*{thm}{Theorem}

\theoremstyle{definition}
\newtheorem{definition}[theorem]{Definition}

\theoremstyle{remark}
\newtheorem{remark}[theorem]{Remark}
\newtheorem{example}[theorem]{Example}

\numberwithin{equation}{section}

\begin{document}
\setlength{\parindent}{0.cm}

\title{String topology coproduct and Turaev cobracket on surfaces}

\author{Jana Hartenstein}
\address{Fachbereich Mathematik, Universität Hamburg, Bundesstraße 55,
20146 Hamburg, Germany}
\email{jana.hartenstein@uni-hamburg.de}

\author{Maximilian Stegemeyer}
\address{Mathematisches Institut, Universit\"at Freiburg, Ernst-Zermelo-Straße 1, 79104 Freiburg, Germany}
\email{maximilian.stegemeyer@math.uni-freiburg.de}
\date{\today}


\begin{abstract}
    The string topology coproduct on the homology of the free loop space of a closed manifold induces a string cobracket on $S^1$-equivariant homology.
    We give a complete computation of the string topology coproduct for surfaces of higher genus by describing an algorithm which computes the coproduct of a cyclic word in terms of generators of the fundamental group of the surface.
    We further show that the string cobracket is the negative of the Turaev cobracket.
\end{abstract}
\maketitle


\section{Introduction}

String topology studies algebraic structures on the homology of the free loop space $LM = \mathrm{Map}(S^1,M)$ of a closed oriented manifold $M$.
In their seminal work \cite{ChasSullivan} Chas and Sullivan introduce  a product on the homology of $LM$ which is known as the \emph{Chas-Sullivan product}.
An operation which can be thought of as a dual to the Chas-Sullivan product is the \emph{Goresky-Hingston coproduct}
$$  \vee \colon H_*(LM,M) \to H_{*+1-n}(LM\times LM, LM\times M\cup M\times LM)   $$
which was first defined and studied in \cite{GoreskyHingston} going back to ideas of Sullivan in \cite{Sul08}.
Here, $M$ is embedded into $LM$ as the constant loops and $n$ is the dimension of $M$.
The Goresky-Hingston coproduct, also known as \emph{string topology coproduct} has received a lot of attention in recent years because of its failure to be homotopy-invariant, see \cite{Naef24}.
Moreover, the question how the Chas-Sullivan product and the string topology coproduct interact with each other turns out to be a delicate problem, see e.g. \cite{CHO20}.

Note that the free loop space is endowed with an $S^1$-action which comes from the reparametrization of loops.
String topology therefore does not only consider algebraic structures on $H_*(LM)$ but also on the equivariant homology $H_*^{S^1}(LM)$.
The Chas-Sullivan product induces a Lie bracket on $H_*^{S^1}(LM)$ while the Goresky-Hingston coproduct induces a cobracket, the \emph{string cobracket} on $H_*^{S^1}(LM)$.
Before Chas and Sullivan's work \cite{ChasSullivan} it had already been known that the equivariant homology of the free loop space of a surface carries a Lie bialgebra structure which consists of the \emph{Goldman bracket} \cite{Goldman} and the \emph{Turaev cobracket} \cite{Turaev}.

In this paper we study the string topology coproduct and cobracket on closed oriented surfaces of positive genus.
In particular, we will prove that the string cobracket is the negative of the Turaev cobracket.
Note that the BV operator and the relationship of the Chas-Sullivan product to the Goldman bracket have been studied for surfaces by Kupers in \cite{Kupers}.

In Section \ref{sec coproduct} we show how one can build a closed oriented surface of genus $g\geq 2$ in such a way that there exists a neighbourhood $U_0$ of the basepoint which is homeomorphic to the $2$-disk and such that for each free homotopy class of loops $[\gamma]$ there is a representative of $[\gamma]$ with all self-intersections only appearing in $U_0$.
This construction allows us to deform the loop so that all self-intersections are transversal and such that we can algorithmically record the self-intersections as well as the associated signs induced by the orientations.
Note that if $c_1,\ldots,c_{2g}\in\pi_1(M)$ are generators of $\pi_1(M)$ then every free homotopy class of loops comes from a cyclic word in the alphabet $\{c_1,\ldots,c_{2g},c_1^-,\ldots,c_{2g}^-\}$ where $c_i^- = c_i^{-1}$ for $i\in \{1,\ldots, 2g\}$.
We show the following.
\begin{thm}[Theorem \ref{theorem_algorithm}]
  Let $M$ be a closed oriented surface of genus $g\geq 2$.
   Let $v = v_1 \ldots v_m$ be a cyclically reduced representative of a cyclic word in the alphabet $\{c_1,\ldots, c_{2g},c_1^-,\ldots ,c_{2g}^- \}$ and let $[\gamma]\in {H}_0(LM)$ be the induced free homotopy class of loops in $M$.
   There is an algorithm taking as input only the letters $v_1,\ldots , v_m$ of $v$ which associates to $v$ finite sets $\mathcal{Q},\overline{\mathcal{Q}}\subseteq \{1,\ldots ,m+1\}^{\times 2}$ as well as signs $\kappa_{j,k}, \overline{\kappa}_{j',k'}\in \{\pm 1\}$ for $(j,k)\in\mathcal{Q}$, $(j',k')\in\overline{\mathcal{Q}}$ such that 
   \begin{eqnarray*}
       \vee \Delta [\gamma] &=& \sum_{(j,k)\in\mathcal{Q}} \kappa_{j,k} \big([v_k\ldots v_{j-1}]\times [v_j\ldots v_{k-1}] - [v_j\ldots v_{k-1}]\times [v_k\ldots v_{j-1}] \big) + \\ & & 
       \sum_{(j,k)\in \overline{\mathcal{Q}}}  \overline{\kappa}_{j,k} \big( [v_{k+1}\ldots v_{j-1}]\times [v_j\ldots v_k] - [v_j\ldots v_k]\times [v_{k+1}\ldots v_{j-1}]\big) .
   \end{eqnarray*}
\end{thm}
Here, $\Delta\colon H_*(LM)\to H_{*+1}(LM)$ is the BV operator and we take the coproduct of the class induced by $\Delta[\gamma]$ in relative homology $H_1(LM,M)$.
The proof of this theorem yields an explicit algorithm which determines the sets $\mathcal{Q},\overline{\mathcal{Q}}$ as well as the signs $\kappa_{j,k}$ and $\overline{\kappa}_{j,k}$.
Further, we will see that the above result yields a full computation of the coproduct on surfaces of genus $g\geq 2$.
In particular, this theorem shows that the coproduct is non-trivial for surfaces of genus $g\geq 2$ as the algorithm yields a non-trivial coproduct for many cyclic words, see also Example \ref{example_coproduct_algorithm}.
However, we show in Theorem \ref{thm repeated simple loop} that the coproduct vanishes for powers of simple loops.
We note that the coproduct for the orientable surfaces of genus $g < 2$ behaves as follows.
The coproduct on the torus with arbitrary coefficients is trivial, see \cite{Stegemeyer}, while the coproduct on the $2$-sphere is non-trivial, see \cite{Ste24} for a computation with rational coefficients.

Let ${L}$ be the set of free homotopy classes of loops in $M$ and let $V$ be the rational vector space spanned by ${L}$.
The Turaev cobracket is a cobracket
$   \vee_T \colon V\to V\otimes V    , $
see \cite{Turaev}.
It turns out that the vector space $V$ is isomorphic to the equivariant homology ${H}_0^{S^1}(LM)$ by sending a representative $\gamma$ of a class in ${L}$ to the corresponding point class $[\gamma]\in {H}_0^{S^1}(LM)$.
We show the following.
\begin{thm}[Theorem \ref{theorem_turaev}]
    Let $M$ be a closed oriented surface of genus $g\geq 2$.
    Under the natural identification $V \cong H_0^{S^1}(LM)$ the string cobracket is the negative of the Turaev cobracket.
\end{thm}
While it seems to have been known to experts that the string cobracket and the Turaev cobracket on closed oriented surfaces are essentially the same, to the best of the authors' knowledge a proof of this statement has not appeared in the literature yet.
Moreover, the above Theorem shows that care has to be taken when it comes to the precise signs coming from the orientation conventions.
We note that we follow the conventions for the string coproduct of Hingston and Wahl in \cite{HingstonWahl}.

\medskip
This article is organized as follows.
In Section 2 we introduce the string topology coproduct as well as the BV operator and recall how the coproduct can be computed geometrically. The homology of the free loop space of surfaces of higher genus is studied in Section 3.
In Section 4 we describe the algorithm which yields a full computation of the string coproduct for surfaces of higher genus.
The comparison between the Turaev cobracket and the string cobracket is carried out in Section 5.
 Finally in Section 6 we give some remarks on the behaviour of the string coproduct on three-manifolds and on hyperbolic manifolds of arbitrary dimension.

\medskip
\noindent \textbf{Acknowledgements.} 
The authors want to thank Nathalie Wahl for many helpful conversations regarding the string topology coproduct. 
In particular, J. H. wants to thank Nathalie Wahl for suggesting the topic.
M. S. was partially funded by the Deutsche Forschungsgemeinschaft (German Research Foundation) -- grant agreement number 518920559.
M. S. is grateful for the support by the Danish National Research Foundation through the Copenhagen Centre for Geometry and Topology (DNRF151).

\section{String topology operations}
\label{sec string operations}

\noindent In this section we define the string topology coproduct following \cite{HingstonWahl} as well as the BV operator. 
We then give a geometric interpretation of the coproduct which will be our tool for computing the coproduct on surfaces.

Throughout this paper we assume that all manifolds are smooth, orientable, closed and connected.
We consider the \emph{free loop space} of a manifold $M$ as the space of all absolutely continuous loops in $M$, i.e. we set
$$  L M = \{\gamma\colon [0,1]\to M\,|\, \gamma(0) = \gamma(1) ,\, \gamma \text{ is of class } H^1 \}.    $$
For the notion of $H^1$-paths we refer to \cite{klingenberg:1995}. 
The space $LM$ can be endowed with the structure of a Hilbert manifold.
Moreover, we note that the space of $H^1$-loops is homotopy equivalent to the space $C^0({S}^1,M)$ of continuous loops in $M$ endowed with the compact-open topology.
For more details, see \cite{HingstonWahl} and \cite{Oancea}.
We shall mostly work with homology with integer coefficients throughout the paper, generalizations to arbitrary coefficients are straightforward.

We fix a Riemannian metric $g$ on $M$ and consider the distance function $\mathrm{d}\colon M\times M\to [0,\infty)$ induced by $g$.
Let $\epsilon >0$ be a number smaller than the injectivity radius of $(M,g)$.
Consider the inclusion of the diagonal $\Delta M\hookrightarrow M\times M$.
A tubular neighborhood is given by 
$$  U_M  = \{ (p,q)\in M\times M\,|\, \mathrm{d}(p,q) < \epsilon \} .  $$
For an $\epsilon_0> 0$ with $\epsilon_0<\epsilon$ we further define the set
$ U_{M,\geq\epsilon_0} = \{ (p,q)\in U_M\,|\, \mathrm{d}(p,q)\geq\epsilon_0\} .  $
Since the inclusion map $(U_M,U_{M,\geq\epsilon_0}) \hookrightarrow (U_M,U_M\setminus \Delta M)$ induces an isomorphism in relative cohomology, the Thom class $\tau\in{H}^n(U_M, U_M\setminus \Delta M) $ induces a class $\tau_M\in{H}^n(U_M,U_{M,\geq\epsilon_0})$.
Consider the evaluation map 
$$
    e_I: LM \times I \xrightarrow[]{} M \times M \qquad 
    e_I(\gamma, s) = (\gamma(0), \gamma(s)) . $$
We set 
$
    U_{\mathrm{GH}}= e_I^{-1}(U_M) $ and $
    \mathcal{F} = e_I^{-1}(\Delta M)$.
Further define $U_{\mathrm{GH},\geq\epsilon_0} = e_I^{-1}(U_{M,\geq\epsilon_0})$.
The evaluation $e_I$ induces a map of pairs $e_I\colon (U_{\mathrm{GH}},U_{\mathrm{GH},\geq\epsilon_0})\to (U_M,U_{M,\geq\epsilon_0})$ and we can thus pull back the class $\tau_M$ to a class $\tau_{\mathrm{GH}}\in {H}^n(U_{\mathrm{GH}},U_{\mathrm{GH},\geq\epsilon_0})$.
Moreover, there is a retraction $R_{\mathrm{GH}}: U_{\mathrm{GH}} \xrightarrow[]{} \mathcal{F}$ which can be defined by using minimizing geodesic sticks, see \cite{HingstonWahl}.
There is further a cutting map $\mathrm{cut}\colon \mathcal{F}\to LM\times LM$ given by cutting the loop $\gamma$ at time $s$ into its two parts for $(\gamma,s)\in\mathcal{F}$.
Finally, let $I=[0,1]$ be the unit interval and let $[I]\in{H}_1(I,\partial I)$ the positively oriented generator represented by the cycle given by $\mathrm{id}_I\colon I\to I$.
\begin{definition}
    Let $M$ be a closed oriented manifold of dimension $n$.
    The \emph{string topology coproduct} is defined as the composition 
    \begin{align*}
        \vee: H_p(LM,M) \xrightarrow[]{\times [I]} H_{p+1}(LM \times I, LM \times \partial I \cup M \times I) \xrightarrow[]{\tau_{GH} \cap} H_{p+1-n}(U_{\mathrm{GH}}, LM \times \partial I \cup M \times I) \\
        \xrightarrow[]{(R_{\mathrm{GH}})_*} H_{p+1-n}(\mathcal{F}, LM \times \partial I \cup M \times I) \xrightarrow[]{\text{cut}_*} H_{p+1-n}(LM \times LM, M \times LM \cup LM \times M).
    \end{align*}
\end{definition}
\begin{remark}
        Hingston and Wahl call the above coproduct the \emph{Thom-signed} coproduct and introduce a sign-corrected version which behaves better when one considers the algebraic properties of the coproduct.
        However for even-dimensional manifolds the sign change is trivial, see \cite[Section 1.5]{HingstonWahl}.
        Since we shall mostly study the coproduct on surfaces in this manuscript, we will not discuss the signs any further.
\end{remark}

\label{subsec geometric coproduct}
Hingston and Wahl provide the following geometric interpretation of the coproduct.

\begin{theorem} [\cite{HingstonWahl}, Prop. 3.7]
\label{theorem_geometric_computation_coproduct}
    Let $(Z: (\Sigma, \Sigma_0) \xrightarrow[]{} (LM,M)) \in C_k(LM, M)$ be a cycle represented by a pair of oriented manifolds $(\Sigma, \Sigma_0)$. Let $\Sigma_B = \Sigma \times \partial I \cup \Sigma_0 \times I$ and 
    \begin{align*}
        E(Z) := e_I \circ Z|_{(\Sigma \times I) \backslash \Sigma_B}:\  (\Sigma \times I)\backslash\Sigma_B & \xrightarrow[]{} M \times M \\
        (\sigma, t) & \mapsto (Z(\sigma)(0), Z(\sigma)(t)).
    \end{align*}
    Assume that $E(Z)$ is tranverse to the diagonal map $\Delta: M \xrightarrow[]{} M \times M$. Let 
    $
        \Sigma_\Delta = E(Z)^{-1}(\Delta M)
    $
    be oriented such that for $(\sigma, t) \in \Sigma_\Delta$ the isomorphism $T_{(\sigma,t)}(\Sigma \times I) \xrightarrow[\cong]{dE(Z)} N_{E(Z)(\sigma,t)} \Delta M \oplus T_{(\sigma,t)}\Sigma_\Delta$
    is orientation-preserving. Then the coproduct of $[Z]$ is given by
    \begin{equation*}
        [\vee Z] = [\mathrm{cut} \circ ((Z \times I)|_{\overline{\Sigma}_\Delta})] \in H_{k+1-n}(LM \times LM, L M \times M \cup M \times L M).
    \end{equation*}
\end{theorem}

Heuristically this theorem says that, in order to compute the coproduct, we only need to consider the loops in the homology class that actually have non-trivial self-intersections, and cut them apart.
Since we will make use of the above theorem later we note that the normal bundle is oriented as follows.
For $(\sigma,t)\in \Sigma_{\Delta}$ with $(p,p) = E(Z)(\sigma,t)\in \Delta  M$ we have an isomorphism
$$    N_{(p,p)} \Delta M  = (T_p M\oplus T_p M) / (\Delta T_p M) \xrightarrow[\cong]{\Phi} T_p M        $$
given by $\Phi[(v,w)] = w-v$ and the orientation of $N_{(p,p)}\Delta M$ is chosen such that this isomorphism are orientation-preserving.

We end this section by introducing the \emph{BV operator} on $H_{*}(LM)$.
There is an $S^1$-action on the free loop space which we denote by $\varphi\colon S^1\times LM\to LM$. 
It is given by
$$
    \varphi(\theta, \gamma)(t) =(\gamma(\theta + t)) \quad \text{for}\,\,\, \theta\in S^1, \,\, \gamma\in LM , \,\,t\in I .
$$
We define the \emph{BV operator}
$
    \Delta: H_*(LM) \xrightarrow[]{} H_{*+1}(LM) $ by $    a \mapsto \varphi_*([S^1] \times a)
$
where $[S^1]\in H_1(S^1)$ is the fundamental class of the circle with its standard orientation.
The BV operator together with the Chas-Sullivan product induces indeed the structure of a BV algebra on $H_*(LM)$. 
For more details we refer to \cite{ChasSullivan} and \cite{cohen2006string}.

\section{The homology of the free loop space of surfaces}
\label{sec homology}

In this section we compute the homology groups of the free loop space of orientable surfaces of genus $g \geq 1$ and determine the action of the BV operator.

Let $M$ be a compact orientable surface of genus $g \geq 1$ with basepoint $p_0\in M$.
For $h \in \pi_1(M)$, let $[h]$ denote the conjugacy class of $h$ in $\pi_1(M)$.
The path components of the free loop space correspond bijectively to the conjugacy classes of the fundamental group of the manifold, see \cite[Theorem 1.6]{Oancea}.
Moreover, the conjugacy classes are in bijection with the set $L$ of free homotopy classes of loops on the manifold. 
We denote the component of $LM$ associated to the conjugacy class $[h]$ for $h\in \pi_1(M)$ by $L_{[h]} M$.
Recall that there is the \emph{free loop fibration} 
\begin{equation}
    \Omega M \xrightarrow{} L M
\xrightarrow[]{\text{ev}_0} M,
\end{equation}
with $\mathrm{ev}_0$ being the evaluation $\mathrm{ev}_0(\gamma ) = \gamma(0)$ and with fiber being the based loop space $\Omega M$. 
By restricting the free loop fibration to the path components of $LM$ and recalling that surfaces of genus $g \geq 1$ are aspherical, one sees that the free loop space $LM$ is homotopy equivalent to a product of Eilenberg-MacLance spaces $L_{[h]}M \simeq K(C_{\pi_1(M)}(h), 1)$, where $C_{\pi_1(M)}(h)$ is the centralizer of $h$ in $\pi_1(M)$ and the product is taken over the conjugacy classes $[h]$ of $\pi_1(M)$.
For more details we refer to \cite{Oancea}, Section 1.4.

Since the fundamental group of the torus $T$ is $\pi_1(T) = \mathbb{Z}^2$, the set of conjugacy classes in $\pi_1(T)$ is $\pi_1(T)$ itself and hence there is a path  component of $LT$ for every $(m,n) \in \mathbb{Z}^2$.
As $\pi_1(T) = \mathbb{Z}^2$ is abelian, we have $C_{\pi_1(T)}(h) = \pi_1(T) = \mathbb{Z}^2$ and thus each path component is a $K(\mathbb{Z}^2,1)$ and thus homotopy equivalent to a torus.
Consequently, we have isomorphisms of groups 
\begin{equation*}
    H_*(LT) \cong \bigoplus_{(m,n) \in \mathbb{Z}^2} H_{*}(T) \cong 
    \begin{cases} 
    \mathbb{Z}[\mathbb{Z}^2] & * = 0,2 \\
    \mathbb{Z}^2[\mathbb{Z}^2] & * = 1 \\
    0 & * \geq 3
    \end{cases}
\end{equation*}

Note that the centralizers of non-trivial elements in hyperbolic groups are virtually cyclic, see \cite[Corollary III.$\Gamma$.3.10]{BH13}. Moreover, there is no torsion in fundamental groups of hyperbolic manifolds, see \cite[Theorem II.4.13]{BH13} and it thus follows that the centralizer of a non-trivial element in a hyperbolic group is infinite cyclic.
If $M$ is a surface of genus $g\geq 2$ we thus get the following.
\begin{lemma}\label{lemma_centralizers_hyperbolic_grps}
    Let $M$ be a surface of genus $g\geq 2$.
    For $h \in \pi_1(M),\ h \neq e$, the centralizer $C_{\pi_1(M)}(h)$ is infinite cyclic and hence there are homotopy equivalences $L_{[h]} M \simeq K(C_{\pi_1(M)} ( h) ) \simeq S^1$.
\end{lemma}
Since the centralizer of the unit element is $\pi_1(M)$ itself, we have $L_{[e]} M \simeq M$.
The homology of $LM$ therefore behaves as follows. 
\begin{equation}
    H_*(LM; \mathbb{Z}) = \bigoplus_{[h] \in \text{Conj}(\pi_1(M))} H_*(L_{[h]} M) \cong 
    \begin{cases}
        \mathbb{Z}[L] & *=0 \\
        H_1(M; \mathbb{Z}) \oplus \mathbb{Z}[\widetilde{L}] &  *=1 \\
        \mathbb{Z} & *=2 \\
        0 & \text{otherwise}
    \end{cases}
\end{equation}
where $\widetilde{L}$ is the set of non-trivial free homotopy classes of loops on $M$.
For the homology relative to the constant loops we thus get
\begin{equation}
    H_*(LM, M; \mathbb{Z}) 
    \cong 
    \begin{cases}
        \mathbb{Z}[\widetilde{L}] & *=0 , 1 \\
        0 & \text{otherwise.}
    \end{cases}
\end{equation}

 Finally, following \cite{Kupers} we study the BV operator on $H_*(LM)$.
 Recall that if $\gamma\in LM$ is a loop then we can take the powers of $\gamma$, i.e. for an integer $k\geq 2$ we define $\gamma^k\in LM$ to be the loop $\gamma^k(t) =\gamma(kt)$ for $t\in [0,1]$.
 We say that $\gamma^k$ is the $k$-th \emph{power} of $\gamma$.
 Let $[h]\in L$ be a non-trivial free homotopy class of loops in $M$ with representative $h\in LM$.
 We define the \emph{level} of $[h]$ to be the largest positive integer $l\in\mathbb{N}$ such that $h$ is freely homotopic to a loop $h'\in LM$ which is the $l$-th power of a loop $\gamma\in LM$.
 We write $l([h]) = l$ for the level of $h$.
Moreover, we define a class $\widetilde{h}\in {H}_1(L_{[h]}M;\mathbb{Z})$ by considering the map $\widetilde{h}\colon {S}^1\to L_{[h]}M$ by  
$$
    \widetilde{h} (\theta)(t) =  h'\Big(t + \frac{\theta}{l([h])}\Big)
$$ 
for $\theta\in S^1, t\in I$. Since $h'$ is an $l([h])$-th power, this is well-defined. 
Kupers argues in \cite{Kupers} that $[\widetilde{h}]:=\widetilde{h}_*([S^1]) \in H_1(L_{[h]}M;\mathbb{Z})$ is a generator where $[S^1]\in H_1(S^1)$ is the fundamental class of the circle.
For $[h]$ the trivial free homotopy class, we define $l([h]) = 0$.
\begin{lemma}\label{lemma_bv_operator}
    Let $M$ be a surface of genus $g\geq 1$.
     Let $[h]\in L$ be a free homotopy class represented by a loop $\gamma$ and consider the homology class $[\gamma]\in H_0(LM)$ induced by $[h]$.
     Then the BV operator satisfies $\Delta[\gamma] = l([h])[\widetilde{h}] $ where $[\widetilde{h}] = \widetilde{h}_*([S^1])$ as above.
 \end{lemma}
 \begin{proof}
    Let $[h]\in L$ be a free homotopy class of level $l([h]) = l\in\mathbb{N}$.
    We represent $[h]$ by an $l$-fold power $h'= \sigma^l$.
    For $m = 1,...,l-1$ define a map $\beta_m \colon S^1\to LM$ by
\begin{equation*}
    \beta_m (\theta) = \Big(t \mapsto h'\big(t + \frac{\theta + m}{l}\big)\Big) \quad \text{for}\,\,\,\theta\in S^1.
\end{equation*}
Since $\beta_m$ is induced by the $S^1$ action, the map $\beta_m$ is homotopic to $\widetilde{h}$ and thus $[\beta_m] = [\tilde{h}]$. 
Therefore, we have 
\begin{equation*}
    \Delta[\gamma] = [\widetilde{h}] + \sum_{m=1}^{l-1} [\beta_m] = l([h]) \, [\widetilde{h}] .
\end{equation*}
It is easy to see that the BV operator is trivial on homology classes coming from the component $L_{[e]}M$ of contractible loops.
This completes the proof.
\end{proof}

\section{Computing the coproduct}
\label{sec coproduct}

In this section, we compute the coproduct for classes of the form $\Delta[\gamma]$ with $[\gamma]\in H_0(LM)$. We shall first prove that $\Delta[\gamma]$ has trivial coproduct for $\gamma$ a power of a simple loop. We then give an algorithm for the coproduct in terms of a cyclically reduced representative of a cyclic word.

Let $M$ be an oriented, closed surface of genus $g \geq 1$. 
We begin with a lemma which allows us to compute the coproduct with Theorem \ref{theorem_geometric_computation_coproduct}.
We use the notation of Theorem \ref{theorem_geometric_computation_coproduct} in the formulation and the proof of the lemma.
\begin{lemma}\label{lemma_coproduct_signs}
    Let $M$ be a closed oriented surface of genus $g\geq 1$ and let $\gamma \in LM$ be a smooth loop with only transverse self-intersections.
    Denote the self-intersection times by $(s_1,t_1),(s_2,t_2),\ldots, (s_m,t_m)\in (0,1)$ with $s_i<t_i$ for $i\in\{1,\ldots, m\}$, i.e. $\gamma(s_i) = \gamma(t_i)$ for $i\in \{1,\ldots,m\}$.
    Consider the map $Z\colon S^1 \to LM$ given by $(Z(s))(t) =  \gamma(s+t)$.
    The induced map $E(Z)\colon S^1\times (0,1) \to M\times M$ is transverse to the diagonal and the coproduct is given by
    $$   \vee [Z] = \sum_{i=1}^m \kappa_i \big(    [\gamma|_{[s_i,t_i]}] \times [\gamma|_{[ t_i, 1+ s_i]}]      - [\gamma|_{[t_i,1+s_i]}] \times [\gamma|_{[s_i,t_i]}] \big)     $$
    where $\kappa_i = -1$ if $(\gamma'(s_i),\gamma'(t_i))$ is positively oriented and $\kappa_i = +1$ else. 
\end{lemma}
\begin{proof}
    We set $\Sigma = S^1$ and we begin by arguing that we can apply Theorem \ref{theorem_geometric_computation_coproduct} to the map $Z\colon \Sigma  \to LM, Z(s)(t) = \gamma(s+t)$.
    It is a direct check that we have $$\Sigma_{\Delta} = E(Z)^{-1}(\Delta M) =  \{ (s_i,t_i-s_i)\,|\, i \in \{1,\ldots, m\}\} \cup \{ (t_i, s_i+1-t_i) \,|\, i \in \{1,\ldots,m\}\}   . $$
    Let $i\in\{1,\ldots, m\}$ and consider the self-intersection point $p = \gamma(s_i)$. We introduce the notation $\tau_i = (s_i,t_i-s_i)\in \Sigma_{\Delta}$ and compute the differential $dE(Z)_{\tau_i}\colon T_{\tau_i}(S^1\times I) \to T_{(p,p)} (M\times M)$ as
    $$    dE(Z)_{\tau_i}(\partial_s) = (\gamma'(s_i),\gamma'(t_i)) \quad \text{and}\quad dE(Z)_{\tau_i}(\partial_t )= (0,\gamma'(t_i)) .        $$
    Here, $\partial_s$ is the global frame on the circle induced by the standard parametrization and $\partial_t$ is the global frame on the unit interval induced by the canonical parametrization.
    Since by assumption all self-intersections are transversal points, we see that $$\mathrm{im}\big( dE(Z)_{\tau_i}\big) \oplus T_{(p,p)}\Delta M = T_{(p,p)} (M\times M) .$$ 
    Consequently, $E(Z)$ is transversal to the diagonal at $\tau_i$.
    Similarly, for $\overline{\tau_i} = (t_i,s_i+1-t_i)$ we see that 
    $$    dE(Z)_{\overline{\tau_i}}(\partial_s) = (\gamma'(t_i),\gamma'(s_i)) \quad \text{and}\quad dE(Z)_{\overline{\tau_i}}(\partial_t )= (0,\gamma'(s_i)) .        $$
    Again this is a complementary subspace to $T_{(p,p)}\Delta M$ and since $i\in\{1,\ldots, m\}$ was arbitrary we see that $E(Z)$ is transverse to the diagonal.
    Now, we consider the projection onto the normal space $\pi\colon T_{(p,p)} (M\times M) \to N_{(p,p)}\Delta M$ as well as the isomorphism $\Phi\colon N_{(p,p)}\Delta M\to T_p M, [(v,w)]\mapsto w-v$.
    By definition $\pi$ and $\Phi$ are orientation-preserving.
    We find that for the composition $f_{\tau_i}= \Phi\circ \pi \circ \big( dE(Z)_{\tau_i}\big) \colon T_{\tau_i}(S^1\times I)\to T_p M$ we have
    $$    f_{\tau_i} (\partial_s)  = \gamma'(t_i)-\gamma'(s_i) \quad \text{and} f_{\tau_i}(\partial_t) = \gamma'(t_i) .    $$
    If the basis $(\gamma'(s_i),\gamma'(t_i))$ is positively oriented then we see that $f_{\tau_i}$ is \emph{orientation-reversing}.
    Hence, $(dE(Z))_{\tau_i}$ is orientation-reversing as well.
    Similarly, one finds that if $(\gamma'(s_i),\gamma'(t_i))$ is positively oriented then $(dE(Z))_{\overline{\tau_i}}$ is \emph{orientation-preserving}.
    In case that the basis $(\gamma'(s_i),\gamma'(t_i))$ is negatively-oriented, then $(dE(Z))_{{\tau_i}}$ is orientation-preserving and $(dE(Z))_{\overline{\tau_i}}$ is orientation-reversing.
    Hence by Theorem \ref{theorem_geometric_computation_coproduct} we obtain a contribution from $\tau_i$ and $\overline{\tau}_i$ of 
    $$  \kappa_i  \big( [\gamma|_{[s_i,t_i]}] \times [\gamma|_{[ t_i, 1+ s_i]}]      - [\gamma|_{[t_i,1+s_i]}] \times [\gamma|_{[s_i,t_i]}]  \big)      $$
    with $\kappa_i$ defined as in the statement of the Lemma.
    Since $i\in \{1,\ldots, m\}$ was arbitrary we obtain the claimed formula for the coproduct.
\end{proof}

It is now a direct consequence that the coproduct of the classes $\Delta [\gamma]$ vanishes when $\gamma$ is a power of a simple loop.
Recall that a loop $S^1 \xrightarrow[]{} M$ is called \emph{simple}, if it is an embedding.
\begin{theorem}
\label{thm repeated simple loop}
    Let $M$ be a closed orientable surface of genus $g\geq 1$.
    Let $\gamma \in LM$ be a simple loop, and let $m \in \mathbb{N}$. Then $\vee (\Delta[\gamma^m]) = 0$.
\end{theorem}

\begin{proof}
    For $m = 1$ the statement follows from \cite[Theorem 3.10]{HingstonWahl}.
    For the rest of the proof we therefore only consider $m\geq 2$.
    Let $\gamma\in LM$ be a simple loop.
    We construct a particular representative of $[\gamma]$.
    Note that the normal bundle of $\gamma$ is an orientable line bundle and hence trivial.
    A tubular neighborhood $U_{\gamma}$ is therefore diffeomorphic to a cylinder $U_{\gamma}\cong{S}^1\times (-1,1)$ with $\gamma$ being embedded as the zero section ${S}^1\times \{0\}\subseteq U_{\gamma}$.
    Consider the $m$-fold power $\gamma^m$ for $m\in\mathbb{N}, m\geq 2$.
    We will now \emph{push} the powers of $\gamma$ in the positive direction inside the cylinder, see Figure 1.
       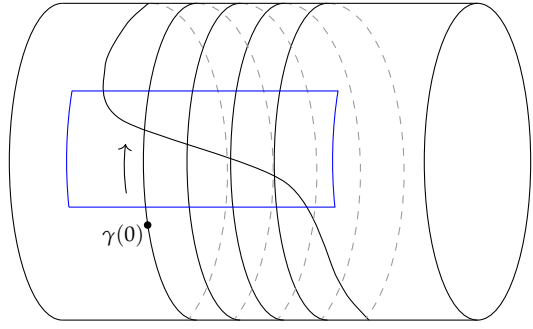
\begin{figure}\label{fig_tube}
	\centering
	\begin{tikzpicture}[scale=1.4]
			\node at (1.3,0.9) [circle,inner sep=1pt,fill=black] {};
			\draw[] (1.76,0) arc (-90:90:-0.5 and 1.5);
		\draw[] (2.17,0) arc (-90:90:-0.5 and 1.5);
		\draw[] (2.58,0) arc (-90:90:-0.5 and 1.5);
		\draw[] (2.99,0) arc (-90:90:-0.5 and 1.5);
		\draw [black] plot [smooth] coordinates { (3.38,0)(3.2,0.2)(3.08,0.4)  (2.6,1.3) (1.05,1.9) (0.9,2.4) (1.1,2.8)  (1.31,3)}; 
		
			\draw[gray!80, dashed] (1.33,3) arc (-90:60:0.72 and -1.6);
			\draw[gray!80, dashed] (1.76,3) arc (-90:60:0.72 and -1.6);
			\draw[gray!80, dashed] (2.17,3) arc (-90:60:0.72 and -1.6);
			\draw[gray!80, dashed] (2.58,3) arc (-90:60:0.72 and -1.6);
			\draw[gray!80, dashed] (2.99,3) arc (-90:60:0.72 and -1.6);
			
			\draw[->] (1.1,1.2) arc (-10:8:-1.3 and 1.5);
			
			 \draw (1.07,0.8) node[scale=0.8]{$\gamma(0)$} ;

        \draw[blue] (0.56,1.07) -- (3.064,1.07);
        \draw[blue] (0.59,2.17) -- (3.09,2.17);
        \draw[blue] (0.56,1.07) arc (-17:26.3:-0.5 and 1.5);
        \draw[blue] (3.06,1.07) arc (-17:26.3:-0.5 and 1.5);
		
		\draw (0.5,0) -- (4.4,0);
		\draw (0.5,3) -- (4.4,3);
		\draw (0.5,0) arc (270:90:0.5 and 1.5);
		\draw (4.4,1.5) ellipse (0.5 and 1.5);
	\end{tikzpicture}
	 \caption{Deformation of $\gamma^m$ in a tubular neighborhood of $\gamma$}
\end{figure}
    More precisely, we let $0<a,b< 1$ with $b > \tfrac{m-1}{m} $ and choose $\epsilon>0$ with $\epsilon < 1-b$.
    Then we define $\widetilde{\gamma}\colon [0,1]\to U_{\gamma}$ to be the loop
    $$   \widetilde{\gamma}(t) = \begin{cases}
        (mt, at), & 0\leq t \leq b \\
        (mt,  \tfrac{a(\epsilon-1)}{\epsilon}t + \tfrac{ab}{\epsilon} ), & b\leq t\leq b+\epsilon  \\
        (mt, at-a), & b+\epsilon \leq t \leq 1.
    \end{cases}   $$
    If we smooth out the corners at $t = b$ and $t= b+\epsilon$ we obtain a smooth loop $\widetilde{\gamma}\in LM$.
    Consider the class $\Delta[\gamma^m]\in H_1(LM,M)$.
    A representative of this class is given by the map $Z\colon S^1 \to LM, \quad (s\mapsto (t\mapsto \widetilde{\gamma}(s+t))  $.
    We want to compute the coproduct of $\Delta [\gamma^m]$ using Lemma \ref{lemma_coproduct_signs}.
    By construction there are times 
    $   0 < s_1 < \ldots < s_{m-1} < t_{m-1} < t_{m-2} < \ldots < t_1 < 1       $
    with $s_i\in (\tfrac{i-1}{m}, \tfrac{i}{m})$ and $t_i\in (b,b+\epsilon)$ for $i\in\{1,\ldots,m-1\}$.
    Let $i\in\{1,\ldots,m-1\}$.
    Cutting the loop $\widetilde{\gamma}_{s_i} := \widetilde{\gamma}(\cdot + s_i) $ at its basepoint yields a pair of loops $(\sigma,\eta)$ with $\sigma$ freely homotopic to $\gamma^{m-i}$ and $\eta$ freely homotopy to $\gamma^i$.
    On the other hand, cutting $\widetilde{\gamma}_{t_i}:= \widetilde{\gamma}(\cdot +t_i)$ at its basepoint yields a pair of loops freely homotopic to the pair $(\gamma^i,\gamma^{m-i})$.
     Hence each pair $(\gamma^i, \gamma^{m-i})$, $i\in\{1,\ldots, m-1\}$ appears exactly twice and by the arguments of Lemma \ref{lemma_coproduct_signs} one sees that the two contributions cancel out.  
\end{proof}

Note that every free homotopy class of loops on the torus can be represented by a power of a simple loop, so the class arising from the BV operator $\Delta$ will have trivial coproduct by Theorem \ref{thm repeated simple loop}.
This observation is already contained in the result by Kupper and the second author who show that the string coproduct on the $r$-torus $T^r$ is trivial for any $r\geq 2$, see Theorem 13 in \cite{Stegemeyer}.

\medskip
We now turn to surfaces of genus $g\geq 2$ and give an algorithm for determining the coproduct.
We will consider free homotopy classes of loops as cyclic words and prove a formula to determine the coproduct in terms of words, where the letters are a set of generators of the fundamental group and their inverses. 

\medskip
\textbf{Surfaces via fattened wedges of circles}

It is well-known that an orientable surface of genus $g\geq 2$ can be obtained by attaching a $2$-disk to a wedge of $2g$ circles.
We shall now construct an orientable surface of genus $g$ by attaching a $2$-disk to a \emph{fattened} wedge of circles.
This idea comes from the study of ribbon graphs and their associated surfaces, see e.g. \cite{Mulase}.

\begin{figure}
	\centering
	\begin{tikzpicture}[scale=1.6]
            \draw[-] (1.5,0) -- (-0.5,0) ;
            \draw[-] (0.5,-1) -- (0.5,1);
             \draw[-] (-0.207,-0.707) -- (1.207,0.707);
             \draw[-] (-0.207,0.707) -- (1.207,-0.707);
   
     \draw (1.7,0.) node[scale=0.8]{$1$} ;
     \draw (1.35,-0.78) node[scale=0.8]{$\overline{4}$} ;
     \draw (0.5,-1.2) node[scale=0.8]{$\overline{3}$} ;
     \draw (-0.35,-0.78) node[scale=0.8]{$\overline{2}$} ;
     \draw (-0.7,0.) node[scale=0.8]{$\overline{1}$} ;
     \draw (-0.35,0.78) node[scale=0.8]{$4$} ;
	\draw (0.5,1.2) node[scale = 0.8]{$3$};
     \draw (1.35,0.78) node[scale=0.8]{$2$} ;
      
      \draw[<-] (0.94,0.2) arc
    [
        start angle=390,
        end angle=150,
        x radius=0.5,
        y radius =0.5
    ] ;

            \draw[-] (5,0) -- (3,0) ;
            \draw[-] (4,-1) -- (4,1);
             \draw[-] (3.293,-0.707) -- (4.707,0.707);
             \draw[-] (3.293,0.707) -- (4.707,-0.707);

             \draw[-, dotted] (4.97,0.1) -- (4.25,0.1);
             \draw[-, dotted] (4.75,0.6) -- (4.25,0.1);
   
            \draw[->] (4.5,0.35) -- (4.6,0.45);
            \draw[->] (4.75,0.1) -- (4.63,0.1);
             
             \draw[-, dotted] (4.97,-0.1) -- (4.25,-0.1);
             \draw[-, dotted] (4.75,-0.6) -- (4.25,-0.1);

             \draw[->] (4.6,-0.45) -- (4.5,-0.35);
                \draw[->] (4.63,-0.1) -- (4.75,-0.1);

             \draw[-, dotted] (3.03,-0.1) -- (3.75,-0.1);
             \draw[-, dotted] (3.25,-0.6) -- (3.75,-0.1);
                
             \draw[->] (3.5,-0.35) -- (3.4,-0.45);
                \draw[->] (3.25,-0.1) -- (3.37,-0.1);
             
             \draw[-, dotted] (3.03,0.1) -- (3.75,0.1);
             \draw[-, dotted] (3.25,0.6) -- (3.75,0.1);
             
             \draw[->] (3.4,0.45) -- (3.5,0.35);
                \draw[->] (3.37,0.1) -- (3.25,0.1);

             \draw[-, dotted] (4.1,0.25) -- (4.1,0.97);
             \draw[-, dotted] (4.1,0.25) -- (4.62,0.75);

             \draw[->] (4.1,0.55) -- (4.1,0.68);
                \draw[->] (4.41,0.55) -- (4.31,0.45);
             
             \draw[-, dotted] (3.9,0.25) -- (3.9,0.97);
             \draw[-, dotted] (3.9,0.25) -- (3.38,0.75);
             
             \draw[->] (3.9,0.68) -- (3.9,0.55);
                \draw[->] (3.69,0.45) -- (3.59,0.55);
             \draw[-, dotted] (4.1,-0.25) -- (4.1,-0.97);
             \draw[-, dotted] (4.1,-0.25) -- (4.62,-0.75);
             
             \draw[->] (3.9,-0.55) -- (3.9,-0.68);
                \draw[->] (3.59,-0.55) -- (3.69,-0.45);
             \draw[-, dotted] (3.9,-0.25) -- (3.9,-0.97);
             \draw[-, dotted] (3.9,-0.25) -- (3.38,-0.75);

             \draw[->] (4.1,-0.68) --  (4.1,-0.55);
                \draw[->]  (4.31,-0.45) --(4.41,-0.55);

              \draw (5.2,0.) node[scale=0.8]{$1$} ;
     \draw (4.85,-0.78) node[scale=0.8]{$\overline{4}$} ;
     \draw (4,-1.2) node[scale=0.8]{$\overline{3}$} ;
     \draw (3.15,-0.78) node[scale=0.8]{$\overline{2}$} ;
     \draw (2.8,0.) node[scale=0.8]{$\overline{1}$} ;
     \draw (3.15,0.78) node[scale=0.8]{$4$} ;
	\draw (4,1.2) node[scale = 0.8]{$3$};
     \draw (4.85,0.78) node[scale=0.8]{$2$} ;

       \draw[dashed] (4,0) circle (0.2cm);
            \draw (4.3,0.12) node[scale=0.8]{$U_0$} ;
	\end{tikzpicture}
        \caption{The wedge of circles $\bigvee_{i=1}^4 {S}^1$ with cyclic order of the half-edges at the vertex and its fattening. The orientation of the boundary of the fattening is indicated. We also sketch the neighborhood $U_0$ of the basepoint.}
        \label{fig_wedge_of_circles_fattening}
\end{figure}
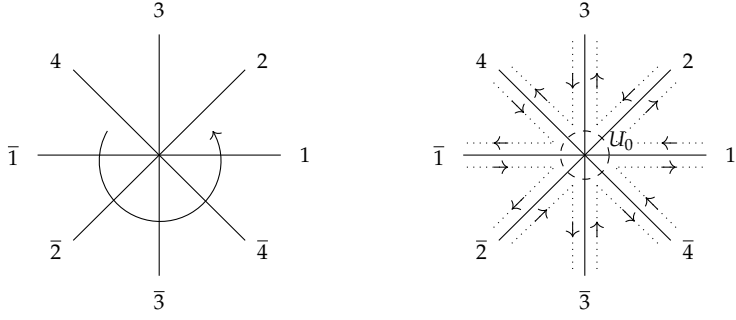

Take a wedge of circles $\bigvee_{i=1}^{2g}S^1$ which we consider as a graph with one single vertex.
We choose a cyclic order on the set of half-edges, see e.g. \cite{ChajdaNovak} for the notion of a cyclic order.
This can be done as follows. Denote the $2g$ circles by $c_1,\ldots, c_{2g}$ and denote the half-edges by $1,\ldots, 2g, \overline{1},\ldots, \overline{2g}$.
On the set $\mathcal{E}$ of half-edges we take the strict total order $1< 2< \ldots 2g<  \overline{1}< \overline{2} < \ldots \overline{2g}$.
A linear order induces a cyclic order, see \cite{ChajdaNovak}, and we take this cyclic order on $\mathcal{E}$, see Figure \ref{fig_wedge_of_circles_fattening}.
We now replace the circles by strips which we connect as indicated in Figure \ref{fig_wedge_of_circles_fattening}.
There is then  precisely one boundary component which we orient by use of the cyclic order, see again Figure \ref{fig_wedge_of_circles_fattening}.
We can smooth out the fattened wedge of circles and obtain an orientable smooth $2$-manifold $\mathcal{F}_{2g}$ with boundary $\partial \mathcal{F}_{2g}\cong S^1$. 
We orient $\mathcal{F}_{2g}$ so that the boundary orientation agrees with the one we have already chosen.
The manifold $\mathcal{F}_{2g}$ is homotopy equivalent to $\bigvee_{i=1}^{2g}S^1$.
We attach a $2$-disk along the boundary using the identification $\partial\mathcal{F}_{2g} \cong S^1$ and the orientation of $\partial\mathcal{F}_{2g}$.
The resulting space $M$ is an orientable closed $2$-manifold, hence it is diffeomorphic to a surface of genus $g'$.
By removing a point in the disk that we attached, we see that the genus $g'$ must indeed be equal to $g$.

As basepoint $p_0\in M$ we choose the basepoint of the wedge of circles.
Note that we have $2g$ distinguished based loops $c_1,\ldots, c_{2g}$ in $M$ given by the circles in the initial wedge of circles.
We shall now study the conjugacy classes in $\pi_1(M)$ and the corresponding homology classes in $LM$.
For the computation of the coproduct we need to control the self-intersections of these loops.
Note that the above construction in particular yields the following.
A loop of the form $\gamma = c_{i_1}^{\pm} \star \ldots \star c_{i_k}^{\pm}$ for $i_1,\ldots, i_k\in\{1,\ldots, 2g\}$ has self-intersections only at the basepoint.
Moreover, by our fattening construction we find a neighbourhood $U_0\subseteq M$ of $p_0$ such that $U_0\subseteq \mathcal{F}_{2g}$ and such that $U_0$ is diffeomorphic to a closed disk of radius $R_0$ in $\mathbb{R}^2$.
We choose this diffeomorphism to be orientation-preserving.
Moreover, the loops $c_1,\ldots, c_{2g}$ intersect $U_0$ in straight lines.
Below we will deform the loop $\gamma$ only in the neighbourhood $U_0$ while leaving it unchanged outside this neighbourhood.
Hence, the only self-intersections that can occur happen in $U_0$ and we can apply Lemma \ref{lemma_coproduct_signs} in order to compute the coproduct of the class $\Delta[\gamma]$.

Next we describe how to understand the free homotopy classes of loops in terms of the generators $c_1,..., c_{2g}$ of $\pi_1(M)$. 
\begin{definition}
\label{def cyclic word}
    Let $A = \{a_1,...,a_q,\overline{a_1},...,\overline{a_q}\}$ be a set of $2q$ letters. A \emph{word} over $A$ is a sequence $v = v_1 \cdots v_l$ with $v_i \in A$ for $i = 1,...,l$. 
    \end{definition}
    We use the convention that $\overline{\overline{a_i}} = a_i$ for $i\in\{1,\ldots q\}$.
    There is an equivalence relation on the words in $A$ generated by 
    \begin{equation*}
        v \sim w \text{ if } w \text{ is a cyclic permutation of } v \text{ or } w = va\overline{a} \text{ for some } a\in A.
    \end{equation*} 
    \begin{definition}
    A \emph{cyclic word} is an equivalence class with respect to this relation. 
    A word $v= v_1\cdots v_l$ is \emph{cyclically reduced} if $v_i \neq \overline{v_{i+1}}$ for all $i = 1,...,l-1$ and $v_l \neq \overline{v_1}$.
\end{definition}

Now we return to the oriented surface of genus $g$ which we considered above and consider the alphabet $A=\{c_1,...,c_{2g},c_1^{-1},...,c_{2g}^{-1}\}$.
Note that for each conjugacy class in $\pi_1(M)$ there is a cyclic word which induces this conjugacy class.
Since the conjugacy classes in $\pi_1(M)$ are in bijective correspondence with the components of $LM$ we obtain a full computation of the coproduct if we compute $\vee \Delta [\gamma]$ for all loops $\gamma$ induced by a cyclic word.

We first need to fix some notation and conventions. 
If $\gamma = v_1 \cdots v_m$ is a cyclically reduced representative of a cyclic word, then we set $v_{m+1 }= v_1$.
Further, if $j>k$ we mean by $v_j\ldots v_k$ the word $v_j\ldots v_m v_1\ldots v_k$.

Recall that we have a fixed identification of the neighborhood $U_0 \subseteq M$ of $p_0\in M$ with the closed disk of radius $R_0$ in $\mathbb{R}^2$.
We choose a smaller radius $0<R_1<R_0$.
For the set $U_0 \setminus \{0\}$ we shall use polar coordinates $(0,R_0]\times S^1\cong U_0\setminus \{0\}$ with $S^1 = \mathbb{R}/\mathbb{Z}$.
For each $i\in\{1,\ldots, 2g\}$ we define points $e_i = (R_0, \tfrac{i-1}{4g})\in U_0$ and $e'_i = (R_0,\tfrac{2g + i-1}{4g})\in U_0$.
We say that $e_i$ and $e_i'$ are the \emph{ends of }$c_i$.
Moreover let $\theta(e_i) = \tfrac{i-1}{4g}$ and $\theta(e'_i) = \tfrac{2g + i-1}{4g}$ for all $i = 1,\ldots, 2g$.
We then have that the intersection of $c_i$ with $U_0$ is the straight line through the origin from $e'_i$ to $e_i$.
We define $\overline{e_i} = e'_i$ and $\overline{e'_i} = e_i$ for all $i\in\{1,\ldots, 2g\}$.

Let $ v_1 \cdots v_m$ be a cyclically reduced representative of a cyclic word in $A$.
By definition there is an $i_j\in \{1,\ldots, 2g\}$ and a $\beta_j\in\{\pm 1\}$ such that $v_j = c_{i_j}^{\beta_j}$ for each $j\in \{1,\ldots ,m\}$.
Hence we write $v_1 \cdots v_m =  c_{i_1}^{\beta_1} \cdots c_{i_m}^{\beta_m}$.
Moreover for each $j\in \{1,\ldots, m\}$ we let $f_j = e_{i_j}$ if $\beta_j = +1$ and $f_j = e'_{i_j}$ else.
Then we can formally write $v_1 \cdots v_m = f_1 \overline{f}_1 f_2 \overline{f}_2 \ldots f_m\overline{f}_m$ and we think of this as the combination as the \emph{outgoing} and the \emph{ingoing} ends of the $c_i$.

Recall that the circle $S^1$ with its standard orientation is a \emph{cyclically orderd set}.
We denote the ternary relation on the circle by $\mathcal{C}_{S^1}\subseteq (S^1)^3$.
Since the set of ends $E=\{e_1,\ldots, e_{2g},e_1',\ldots , e_{2g}'\}$ is defined as a subset of the circle of radius $R_0$ in $U_0$, the set $E$ inherits a cyclic order from the circle.
We denote the ternary relation by $\mathcal{C}_{4g}\subseteq E^3$.
Note that the cyclically ordered set $E$ is isomorphic as a cyclically ordered set to the cyclically ordered set of half-edges $\mathcal{E}$ that we introduced earlier.

\begin{theorem}\label{theorem_algorithm}
   Let $M$ be a closed oriented surface of genus $g\geq 2$.
   Let $v = v_1 \ldots v_m$ be a cyclically reduced representative of a cyclic word in the alphabet $\{c_1,\ldots, c_{2g},c_1^-,\ldots ,c_{2g}^- \}$ and let $[\gamma]\in {H}_0(LM)$ be the induced free homotopy class of loops in $M$.
   There is an algorithm taking as input only the letters $v_1,\ldots , v_m$ which associates to $v$ finite sets $\mathcal{Q},\overline{\mathcal{Q}}\subseteq \{1,\ldots ,m+1\}^{\times 2}$ as well as signs $\kappa_{j,k}, \overline{\kappa}_{j',k'}\in \{\pm 1\}$ for $(j,k)\in\mathcal{Q}$, $(j',k')\in\overline{\mathcal{Q}}$ such that
   \begin{eqnarray*}
       \vee \Delta [\gamma] &=& \sum_{(j,k)\in\mathcal{Q}} \kappa_{j,k} \big([v_k\ldots v_{j-1}]\times [v_j\ldots v_{k-1}] - [v_j\ldots v_{k-1}]\times [v_k\ldots v_{j-1}] \big) + \\ & & 
       \sum_{(j,k)\in \overline{\mathcal{Q}}}  \overline{\kappa}_{j,k} \big( [v_{k+1}\ldots v_{j-1}]\times [v_j\ldots v_k] - [v_j\ldots v_k]\times [v_{k+1}\ldots v_{j-1}]\big) .
   \end{eqnarray*}
\end{theorem}
Note that $\Delta[\gamma]$ is a homology class in $H_1(LM)$, while the coproduct is defined on relative homology $H_1(LM,M)$.
In the statement of the theorem we mean by $\Delta[\gamma]$ the class induced by $\Delta[\gamma]$ in $H_1(LM,M)$ through the natural map $H_1(LM)\to H_1(LM,M)$.
Similarly, the output of the coproduct $\vee\Delta[\gamma]$ lies in $H_0(LM,M)\otimes H_0(LM,M)$ and thus by a class of the form $[v_k\ldots v_{j-1}]$ we mean the class induced by the free homotopy class of $v_k\ldots v_{j-1}$ in $H_0(LM,M)$.
\begin{proof}
Let $v = v_1\ldots v_m$ be a cyclically reduced representative of a cyclic word.
We shall construct a representative of $v$ which has only transverse self-intersections all of which appear in the set $U_0$.
In the course of this proof we describe an algorithm that records all the self-intersections of this particular representative.
In order to compute the coproduct using Lemma \ref{lemma_coproduct_signs} we also need to take care of the respective signs at the self-intersections.
    
As before let $R_1< R_0$.
We write $v = c_{i_1}^{\kappa_1} \ldots c_{i_m}^{\kappa_m}$ where $i_j\in \{1,\ldots, 2g\}$ for $j\in \{1,\ldots, m\}$ and $\kappa_j\in \{\pm 1\}$ for $j\in\{1,\ldots, m\}$.
We write each $c_{i_j}^{\kappa_j}$ as a combination of an outgoing end and an ingoing end $f_{j}\overline{f}_{j}$ as explained above.
We further choose radii $R_1 < r_2 < r_3 <\ldots  < r_m < R_0$ and we choose numbers $0< \epsilon_1 < \ldots  < \epsilon_m < \tfrac{1}{8g}$.
We now describe how we deform the loop $ c_{i_1}^{\kappa_1} \star \ldots \star c_{i_m}^{\kappa_m}$ while keeping track of signs and intersection points.

\medskip
\textbf{Step $\mathbf{1}$:} Define $\delta_1\colon I\to U_0$ to be the path given by
$  
\delta_1(t) = ( t\cdot R_0, \theta(f_{1}) + \epsilon_1) .
$
Define the sets $\mathcal{Q}_1 = \emptyset$ and $\overline{\mathcal{Q}}_1 = \emptyset$.

\medskip
\textbf{Steps $\mathbf{2,\ldots, m}$:}
Let $j\geq 2$.
We want to define a path from the point $(R_0, \theta(\overline{f}_{{j-1}})-\epsilon_j)$ to $(R_0,\theta(f_{j}) + \epsilon_j)$ by first going on a straight line in direction of the origin until we reach radius $r_j$, then going along a circle segment of radius $r_j$ to the angle $\theta(f_{j}) + \epsilon_j$ and then going on a radial line to the endpoint.
We need to choose whether the circle segment goes in the mathematically positive or negative way.
We make the following choice: If $f_{j} \neq f_{{j-1}}$, then we go in the direction of the \emph{shorter way} between $\overline{f}_{j-1}$ and $f_j$, otherwise we go in mathematically positive direction.
Note that by our choice of a cyclically reduced representative we cannot have $f_{j} = \overline{f}_{{j-1}}$.
If $f_{j}\neq f_{{j-1}}$, let $\xi_j\in \mathbb{R}$ be such that the path
\begin{equation}\label{eq_choice_of_path_on_circle}  \Xi_j\colon I\to \mathbb{R} , \quad t\mapsto  (1-t)\cdot (\theta(f_{{j-1}}) - \epsilon_{j-1})   + t \cdot \xi_j      \end{equation}
describes the shortest path on the circle from $\theta(\overline{f}_{j-1})-\epsilon_{j-1}$ to $\theta(f_j) +\epsilon_j$ if we mod out by $\mathbb{Z}$.
Note that in general we cannot choose $\xi_j  = \theta(f_j)+\epsilon_j$ since this convex combination in $\mathbb{R}$ might not yield the shortest way on the circle $\mathbb{R}/\mathbb{Z}$.
If $\Xi_j$ goes in mathematically positive direction we define $s_j = +1$, otherwise we set $s_j = -1$.
If $f_{j} = f_{j-1}$, then we define $\xi_j = \theta(\overline{f}_{j-1}) + \tfrac{1}{2} + (\epsilon_j-\epsilon_{j-1})$ and we define $\Xi_j\colon I\to \mathbb{R}$ as in equation \eqref{eq_choice_of_path_on_circle}.
In this case, equation $\eqref{eq_choice_of_path_on_circle}$ describes a path which goes in positive direcition on the circle.
We further set $s_j = +1$ in this case.
We define the path $\mu_j\colon I\to U_0$
$$   \mu_j(t) = \begin{cases}
    (     3t \cdot r_j +   (1-3t)\cdot R_0 ,     \theta(\overline{f}_{{j-1}}) - \epsilon_{j-1}) , &  0\leq t \leq \tfrac{1}{3} \\
    (r_j, \Xi_j(3t-1)) , & \tfrac{1}{3} \leq t \leq \tfrac{2}{3} \\
    (   (3t-2)\cdot R_0 + (3-3t)\cdot r_2,      \theta(f_{j}) + \epsilon_j), & \tfrac{2}{3} \leq t \leq 1 .
\end{cases}    $$
For a sketch of the path $\mu_j$ we refer to Figure \ref{fig_exampleA}.
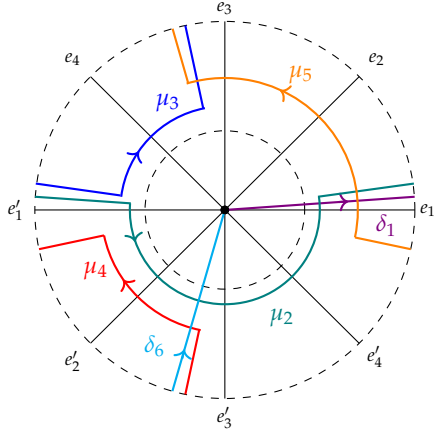
\begin{figure}
	\centering
	\begin{tikzpicture}[scale=2.5]
             \filldraw (0.,0) circle (0.02cm);
            \draw[dashed] (0,0) circle (1.cm);
            \draw[dashed] (0.,0) circle (0.42cm);

            \draw[-] (1.,0) -- (-1,0) ;
           \draw[-] (-0.707,-0.707) -- (0.707,0.707);
            \draw[-] (0.,-1) -- (0,1);
            \draw[-] (-0.707,0.707) -- (0.707,-0.707);

            \draw (1.07,0) node[scale=0.7]{$e_1$};
            \draw (0.8,0.8) node[scale = 0.7]{$e_2$};
            \draw (0.,1.07) node[scale = 0.7]{$e_3$};
            \draw (-0.8,0.8) node[scale = 0.7]{$e_4$};
            \draw (-1.1,0) node[scale=0.7]{$e_1'$};
            \draw (-0.8,-0.8) node[scale = 0.7]{$e_2'$};
            \draw (-0.,-1.1) node[scale = 0.7]{$e_3'$};
            \draw (0.79,-0.777) node[scale = 0.7]{$e_4'$};

            \draw[-, thick ,violet] (0,0) -- (0.999,0.07);
             \draw [->, thick ,violet] (0.6,0.042) -- (0.65,0.0445);
     \draw[violet] (0.85,-0.08) node[scale=0.89]{$\delta_1$} ;

            \draw[-, thick ,teal] (-0.999,0.07) -- (-0.495,0.035)    ; 
             \draw[-, thick, teal] (-184:0.5) arc (-184:8:0.5) ;
            \draw[-,thick, teal] (0.49,0.07) -- (0.995,0.14);
             \draw [->, thick, teal]  (-0.474,-0.16) --  (-0.472,-0.17);
        \draw[teal] (0.3,-0.57) node[scale=0.89]{$\mu_2$} ;


        \draw[-, thick, blue] (-0.995,0.14) -- (-0.55, 0.077);
         \draw[-, thick, blue] (101:0.55) arc (101:172.7:0.55) ;
         \draw[-, thick, blue] (-0.208,0.978) -- (-0.114,0.54);
             \draw [->, thick, blue]  (-0.459,0.3) --  (-0.452,0.31);
            \draw[blue] (-0.3,0.57) node[scale=0.89]{$\mu_3$} ;


        \draw[-, thick, red] (-0.208,-0.978) -- (-0.135,-0.636);
        \draw[-, thick, red] (192:0.65) arc (192:258.7:0.65);
        \draw[-, thick, red] (-0.978,-0.208) -- ( -0.636,-0.135);
         \draw [->, thick, red]  (-0.528,-0.38) --  (-0.536,-0.37);
            \draw[red] (-0.68,-0.32) node[scale=0.89]{$\mu_4$} ;


        \draw[-, thick ,orange] (0.978,-0.208) -- (0.685, -0.146);
        \draw[-, thick, orange] (-12:0.7) arc (-12: 106.7:0.7);
        \draw[-, thick ,orange] (-0.275, 0.961) -- (-0.193,0.673);
         \draw [->, thick, orange]  (0.3,0.63) --  (0.29,0.635);
            \draw[orange] (0.4,0.72) node[scale=0.89]{$\mu_5$} ;

            \draw[-, thick, cyan] (-0.276,-0.961) -- (0,0);
            \draw[->, thick, cyan] (-0.224, -0.788) -- ( -0.213 ,-0.739);
            \draw[cyan] (-0.37,-0.72) node[scale=0.89]{$\delta_6$} ;
        
        \filldraw (0.,0) circle (0.02cm);
        
	\end{tikzpicture}

        \caption{Construction of the loop respresenting the word $c_1c_1c_3c_1^{-1} c_3$ for genus $g= 2$. The following cases of intersections appear: $\mu_2$ intersects $\delta_1$, this is case 1.c).
        Further, $\mu_5$ intersects $\delta_1$ and $\mu_2$, both are case 1.e) and $\mu_5$ intersects with $\mu_3$, this is case 1.c). Finally, $\delta_6$ intersects $\mu_2$, this is case 3.a) and $\delta_6$ and $\mu_4$ intersect giving an example of case 3.d).
         }
        
        \label{fig_exampleA}
    \end{figure}
We now record the intersections of $\mu_j$ with the paths $\delta_1,\mu_2,\ldots, \mu_{j-1}$, i.e. these are the \emph{new} self-intersections occuring at step $j$.
We refer to Figure \ref{fig_exampleA} as well for a sketch where some of the cases that we describe next appear.

By construction the intersections of $\mu_j$ with the previously defined paths can happen only for $\tfrac{1}{3} < t <\tfrac{2}{3}$.
Intersections can only come from the intersection of the segment $(r_j,\Xi_j)$ with the radial lines of the previous paths.
More precisely there is an intersection exactly if there is a radial line with angle $\alpha$ such that the following holds.
\begin{enumerate}
    \item[1.] The orientation of $\Xi_j$ is positive, i.e. $s_j = +1$ and $(\theta(\overline{f}_{j-1}) - \epsilon_{j-1}, \alpha, \theta(f_j) + \epsilon_j)\in \mathcal{C}_{{S}^1}$ or

    \item[2.] The orientation of $\Xi_j$ is negative, i.e. $s_j = -1$ and $(\theta(f_{j}) + \epsilon_j, \alpha, \theta(\overline{f}_{{j-1}} )   - \epsilon_{j-1})\in \mathcal{C}_{{S}^1}$.
\end{enumerate}
We recall that $\mathcal{C}_{S^1}\subseteq (S^1)^3$ is the cyclic order.
By considering the construction of $\delta_1$ and the $\mu_k$, $k< j$ we see that in case 1. we have intersections precisely in the following cases
\begin{enumerate}
    \item[1.a)] there is a $k< j $ such that $(\overline{f}_{j-1},f_{k},f_{j})\in \mathcal{C}_{4g}$ or
    \item[1.b)] there is a $k< j $ such that $(\overline{f}_{j-1},\overline{f}_{k},f_{j})\in \mathcal{C}_{4g}$ or
    \item[1.c)] there is a $k< j$ such that $f_{k} = f_{ j}$ or
    \item [1.d)] there is a $k< j$ such that $\overline{f}_k = f_j$
    \item[1.e)]  there is a $k< j-1$ such that $f_k = \overline{f}_{j-1}$
    \item[1.f)]there is a $k<j-1$ such that $\overline{f}_{k} = \overline{f}_{{j-1}}$
\end{enumerate}
Note that the pair $(j,k)$ can only satisfy at most one of the cases 1.a),1.c) and 1.e) and it can only satisfy at most one of the cases 1.b), 1.d) and 1.f).
We define the sets $$\mathcal{S}_j = \{(j,k)\,|\, k \text{ satisfies condition 1.a) or 1.c) or 1.e)\}} $$ and $$ \overline{\mathcal{S}}_j = \{(j,k) \,|\, k \text{ satisfies condition 1.b) or 1.d) or 1.f)}\}  .   $$
We now need to record the signs, resp. the orientation of the basis vectors at the intersection points.
This can be done case by case, see Figure \ref{fig_self_intersection}.
If $(j,k)\in\mathcal{S}_j$  we have a positively oriented basis at the intersection points and hence we define $\kappa_{j,k} = -1$.
If $(j,k)\in\overline{\mathcal{S}}_j$ we have a negatively oriented basis and therefore we set $\overline{\kappa}_{j,k} = +1$.
Define $\mathcal{Q}_j = \mathcal{Q}_{j-1} \cup \mathcal{S}_j$ and $\overline{\mathcal{Q}}_j = \overline{\mathcal{Q}}_{j-1} \cup \overline{\mathcal{S}}_j $.

In case 2. we have intersections if and only if one of the following holds
\begin{enumerate}
    \item[2.a)] there is a $k< j$ such that $(f_{j},f_{k},\overline{f}_{{j-1}})\in \mathcal{C}_{4g}$ or
    \item[2.b)] there is a $k< j $ such that $(f_{j}, \overline{f}_{k},\overline{f}_{{j-1}})\in \mathcal{C}_{4g}$.
\end{enumerate}
We define the set
$$    \mathcal{S}_j =  \{ (j,k)\,|\, k \text{ satisfies condition 2.a)}\}      \quad \text{and}\quad  \overline{\mathcal{S}}_j = \{ (j,k)\,|\, k\text{ satisfies condition 2.b)}\}   . $$
Furthermore we set $\kappa_{j,k} = +1$ if $(j,k)\in\mathcal{S}_j$ and $\overline{\kappa}_{j,k} = -1$ if $(j,k)\in\overline{\mathcal{S}}_j$.
Define $\mathcal{Q}_j = \mathcal{Q}_{j-1} \cup \mathcal{S}_j$ and $\overline{\mathcal{Q}}_j = \overline{\mathcal{Q}}_{j-1} \cup \overline{\mathcal{S}}_j $.

\medskip
\textbf{Step $\mathbf{m+1}$:}
Define the line $\delta_{m+1}\colon I\to U_0$ by 
$   \delta_{m+1}(t) = ( (1-t)\cdot R_0, \theta(\overline{f}_{m}) - \epsilon_{m})    $.
We now obtain self-intersections if there is a $k\in \{2,\ldots, m\}$ such that $(r_k,\Xi_k)$ and $\delta_{m+1}$ intersect.
One checks that this happens precisely in one of the following four cases:
\begin{enumerate}
    \item[3.a)] there is a $k\in\{2,\ldots, m-1\}$ such that $s_k = +1$ and $(\overline{f}_{{k-1}}, \overline{f}_{m}, f_{k}) \in \mathcal{C}_{4g}$ or
    \item[3.b)] there is a $k\in \{2,\ldots, m-1\} $ such that $s_k = -1$ and $(f_{k},\overline{f}_{m},\overline{f}_{{k-1}})\in\mathcal{C}_{4g}$ or
     \item[3.c)] there is a $k\in\{2,\ldots, m-1\}$ such that $s_k = +1$ and $f_{k} = \overline{f}_{m}$ or
     \item[3.d)] there is a $k\in\{2,\ldots, m-1\}$ such that $s_{k} = -1$ and $\overline{f}_{{k-1}} = \overline{f}_{m}$. 
\end{enumerate}
We define the set
$$  \mathcal{S}_{m+1} = \{ (m+1,k) \,|\, k \text{ satisfies one the conditions 3.a) - 3.d)}\} .   $$ 
For $(m+1,k)\in\mathcal{S}_{m+1}$ we set $\kappa_{m+1,k} = -1$ if $k$ satisfies condition 3.a) or 3.c) and $\kappa_{m+1,k} = +1$ if $k$ satisfies condition 3.b) or 3.d).
Finally, define $\mathcal{Q} := \mathcal{Q}_m\cup \mathcal{S}_{m+1}$ and $\overline{\mathcal{Q}}:= \overline{\mathcal{Q}}_m$.

\medskip
We now define a loop $\gamma\in LM$ representing the free homotopy class given by the cyclically reduced word $v$.
For $i\in\{1,\ldots,m\}$ define a path $\nu_j\colon I\to M$ which connects the points $(R_0,\theta(f_{j})+ \epsilon_j)$ and $(R_0,\theta(\overline{f}_{j})-\epsilon_j)$ by going along the strip which we use in the construction of the surface.
We then define $\gamma$ as the concatenation
$$    \gamma = \delta_1 \star \nu_1 \star \mu_2 \star \nu_2 \star \mu_3 \star \cdots \star \mu_{m} \star \nu_m \star \delta_{m+1} .     $$
By construction there are no intersections coming from the paths $\nu$ so we have captured all self-intersections of $\gamma$ in the above discussion.
The corners of $\gamma$ can be smoothened out without creating new self-intersections.
Moreover by construction all self-intersections are transversal.
The computation of the coproduct follows now directly from Lemma \ref{lemma_coproduct_signs} and the fact that we defined the signs $\kappa_{j,k}$ following Lemma \ref{lemma_coproduct_signs}.

\begin{figure}
	\centering
	\begin{tikzpicture}[scale=1.6]       
   
     \draw[-] (0.0,-0.9) -- (0,-1.7);
     \draw [->] (0.,-1.29) -- (0.,-1.3);

    \draw[dashed] (0.75,-1.3) -- (1.0,-1.56);
    \draw[dashed] (-0.75,-1.3) -- (-1,-1.56);
     \draw[-] (240:1.5) arc (240:300:1.5) ;
     \draw [->] (0.4,-1.45) -- (0.42,-1.44);

     \draw (0.6,-1.18) node[scale=0.7]{$\mu_j$} ;
     \draw (-0.1,-1.05) node[scale=0.7]{$f_k$} ;

     \draw (-0.3,-1.9) node[scale=1.]{\textbf{1.a)}, $\kappa_{j,k} = -1$} ;


    \draw[-] (2.8,-0.9) -- (2.8,-1.7);
     \draw [->]  (2.8,-1.3) -- (2.8,-1.29);

    \draw[dashed] (3.55,-1.3) -- (3.8,-1.56);
    \draw[dashed] (2.05,-1.3) -- (1.8,-1.56);
     \draw[-](2.05,-1.3) arc (240:300:1.5) ;
     \draw [->]  (3.2,-1.45) -- (3.22,-1.44) ;

     \draw (3.4,-1.18) node[scale=0.7]{$\mu_j$} ;
     \draw (2.7,-1.05) node[scale=0.7]{$\overline{f}_k$} ;

     \draw (2.5,-1.9) node[scale=1.]{\textbf{1.b)}, $\overline{\kappa}_{j,k} = +1$} ;


    \draw[dotted] (5.9,-0.9) -- (6.05,-1.8);
    \draw (5.8,-1.) node[scale=0.7]{$e_*$};
     
    \draw[] (6.03,-1.) -- (6.25,-1.8);
     \draw [->] (6.108,-1.29)  -- (6.111,-1.3);

    \draw[] (6.35,-1.3) -- (6.6,-1.56);
    \draw[dashed] (4.85,-1.3) -- (4.6,-1.56);
     \draw[-](4.85,-1.3) arc (240:300:1.5) ;
     \draw [->]  (6.27,-1.34) -- (6.29,-1.33) ;

     \draw (6.6,-1.38) node[scale=0.7]{$f_j$} ;
     \draw (6.17,-1.05) node[scale=0.7]{${f}_k$} ;

     \draw (5.3,-1.9) node[scale=1.]{\textbf{1.c)}, ${\kappa}_{j,k} = -1$} ;


     \draw[] (0.325,-2.55) -- (0.45,-3.3);
      \draw [->] (0.363,-2.79)  -- (0.361,-2.78);
    \draw (0.2,-2.5) node[scale=0.7]{$\overline{f}_k$};
     
    \draw[dotted] (0.43,-2.5) -- (0.65,-3.3);
   
    \draw[] (0.8,-2.8) -- (1,-3.06);
    \draw[dashed] (-0.75,-2.8) -- (-1,-3.06);
     \draw[-](-0.7,-2.8) arc (240:300:1.5) ;
     \draw [->]  (0.67,-2.87) -- (0.69,-2.86) ;

     \draw (1,-2.88) node[scale=0.7]{$f_j$} ;
     \draw (0.57,-2.55) node[scale=0.7]{$e_*$} ;

     \draw (-0.3,-3.4) node[scale=1.]{\textbf{1.d)}, ${\overline{\kappa}}_{j,k} = +1$} ;


       \draw[] (2.475,-2.55) -- (2.375,-3.2);
      \draw [->] (2.439,-2.78)  -- (2.437,-2.79);
    \draw (2.6,-2.5) node[scale=0.7]{${f}_k$};
     
    \draw[dotted] (2.37,-2.5) -- (2.175,-3.2);
   
    \draw[] (2.1,-2.8) -- (1.9,-3.06);
    \draw[dashed] (3.6,-2.8) -- (3.85,-3.06);
     \draw[-](2.1,-2.8) arc (240:300:1.5) ;
     \draw [->]  (2.61,-2.985) -- (2.63,-2.99) ;

     \draw (1.8,-2.88) node[scale=0.7]{$\overline{f}_{j-1}$} ;
     \draw (2.23,-2.55) node[scale=0.7]{$e_*$} ;

     \draw (2.5,-3.4) node[scale=1.]{\textbf{1.e)}, ${{\kappa}}_{j,k} = -1$} ;


       \draw[dotted] (5.275,-2.55) -- (5.175,-3.2);
      \draw (5.4,-2.5) node[scale=0.7]{$e_*$};
     
    \draw[] (5.17,-2.5) -- (4.975,-3.2);
   \draw [->]  (5.1,-2.74) -- (5.102,-2.73);
    
    \draw[] (4.9,-2.8) -- (4.7,-3.06);
    \draw[dashed] (6.4,-2.8) -- (6.65,-3.06);
     \draw[-](4.9,-2.8) arc (240:300:1.5) ;
     \draw [->]  (5.41,-2.985) -- (5.43,-2.99) ;

     \draw (4.6,-2.88) node[scale=0.7]{$\overline{f}_{j-1}$} ;
     \draw (5.03,-2.55) node[scale=0.7]{$\overline{f}_k$} ;

     \draw (5.3,-3.4) node[scale=1.]{\textbf{1.f)}, ${\overline{\kappa}}_{j,k} = +1$} ;


             \draw[-] (0.0,-3.9) -- (0,-4.7);
     \draw [->] (0.,-4.29) -- (0.,-4.3);

    \draw[dashed] (0.75,-4.3) -- (1.0,-4.56);
    \draw[dashed] (-0.75,-4.3) -- (-1,-4.56);
     \draw[-](-0.75,-4.3) arc (240:300:1.5) ;
     \draw [->] (-0.4,-4.45) -- (-0.42,-4.44);

     \draw (-0.6,-4.18) node[scale=0.7]{$\mu_j$} ;
     \draw (-0.1,-4.05) node[scale=0.7]{$f_k$} ;

     \draw (-0.3,-4.9) node[scale=1.]{\textbf{2.a)}, $\kappa_{j,k} = +1$} ;


     \draw[-] (2.8,-3.9) -- (2.8,-4.7);
     \draw [->]  (2.8,-4.3) -- (2.8,-4.29);

    \draw[dashed] (3.55,-4.3) -- (3.8,-4.56);
    \draw[dashed] (2.05,-4.3) -- (1.8,-4.56);
     \draw[-](2.05,-4.3) arc (240:300:1.5) ;
     \draw [->]  (2.4,-4.45) -- (2.38,-4.44) ;

     \draw (3.4,-4.18) node[scale=0.7]{$\mu_j$} ;
     \draw (2.7,-4.05) node[scale=0.7]{$\overline{f}_k$} ;

     \draw (2.5,-4.9) node[scale=1.]{\textbf{2.b)}, $\overline{\kappa}_{j,k} = -1$} ;


       \draw[-] (5.6,-3.9) -- (5.6,-4.7);
     \draw [->]  (5.6,-4.3) -- (5.6,-4.29);

    \draw[dashed] (6.35,-4.3) -- (6.6,-4.56);
    \draw[dashed] (4.85,-4.3) -- (4.6,-4.56);
     \draw[-](4.85,-4.3) arc (240:300:1.5) ;
     \draw [->]  (6.,-4.45) -- (6.02,-4.44) ;

     \draw (6.2,-4.18) node[scale=0.7]{$\mu_k$} ;
     \draw (5.45,-4.05) node[scale=0.7]{$\overline{f}_m$} ;

     \draw (5.3,-4.9) node[scale=1.]{\textbf{3.a)}, ${\kappa}_{m+1,k} = -1$} ;


     \draw[-] (0.0,-5.4) -- (0,-6.2);
     \draw [->]  (0.,-5.8)  -- (0.,-5.79);

    \draw[dashed] (0.75,-5.8) -- (1.0,-6.06);
    \draw[dashed] (-0.75,-5.8) -- (-1,-6.06);
     \draw[-] (-0.75,-5.8) arc (240:300:1.5) ;
     \draw [->] (-0.4,-5.95) -- (-0.42,-5.94);

     \draw (0.6,-5.68) node[scale=0.7]{$\mu_k$} ;
     \draw (-0.15,-5.55) node[scale=0.7]{$\overline{f}_m$} ;

     \draw (-0.3,-6.4) node[scale=1.]{\textbf{3.b)}, $\kappa_{m+1,k} = + 1$} ;


     \draw[] (3.125,-5.55) -- (3.225,-6.2);
      \draw [->] (3.163,-5.79)  -- (3.161,-5.78);
    \draw (2.95,-5.5) node[scale=0.7]{$\overline{f}_m$};
     
    \draw[dotted] (3.23,-5.5) -- (3.45,-6.3);
   
    \draw[] (3.6,-5.8) -- (3.8,-6.06);
    \draw[dashed] (2.05,-5.8) -- (1.8,-6.06);
     \draw[-](2.1,-5.8) arc (240:300:1.5) ;
     \draw [->]  (3.47,-5.87) -- (3.49,-5.86) ;

     \draw (3.85,-5.88) node[scale=0.7]{$f_k$} ;
     \draw (3.37,-5.55) node[scale=0.7]{$e_*$} ;

     \draw (2.5,-6.4) node[scale=1.]{\textbf{3.c)}, ${{\kappa}}_{m+1,k} = -1$} ;


     \draw[dotted] (6.1,-5.4) -- (6.5,-6.3);
    \draw (6.25,-5.5) node[scale=0.7]{$e_*$};
     
    \draw[] (5.85,-5.5) -- (5.95,-6.2);
     \draw [->] (5.89,-5.8) --  (5.888,-5.79);

    \draw[] (6.13,-5.9) -- (6.25,-6.16);
    \draw[dashed] (4.85,-5.8) -- (4.6,-6.06);
     \draw[-](4.85,-5.8) arc (240:290.7:1.5) ;
     \draw [->] (5.7,-6) --  (5.68,-6) ;

    \draw (6.1,-6.08) node[scale=0.7]{$f_k$} ;
     \draw (5.72,-5.55) node[scale=0.7]{$\overline{f}_m$} ;

     \draw (5.3,-6.4) node[scale=1.]{\textbf{3.d)}, ${\kappa}_{m+1,k} = +1$} ;
     
	\end{tikzpicture}
        \caption{Sketch of the twelve cases of intersections appearing in the algorithm. For each case we record the corresponding sign which is determined by Lemma \ref{lemma_coproduct_signs}. Note that in the cases 1.a) - 1.f) as well as 2.a) and 2.b) the part $f_k$ or $\overline{f}_k$ is run through first by the loop before $\mu_j$. In the cases 3.a) - 3.d) the part $\mu_k$ is run through first before $\overline{f}_m$. For cases 1.c) - 1.f) as well as 3.c) and 3.d) where the ends agree, we draw the original position of the end $e_*$ in order to point out where the intersections happen.}
        \label{fig_self_intersection}
    \end{figure}
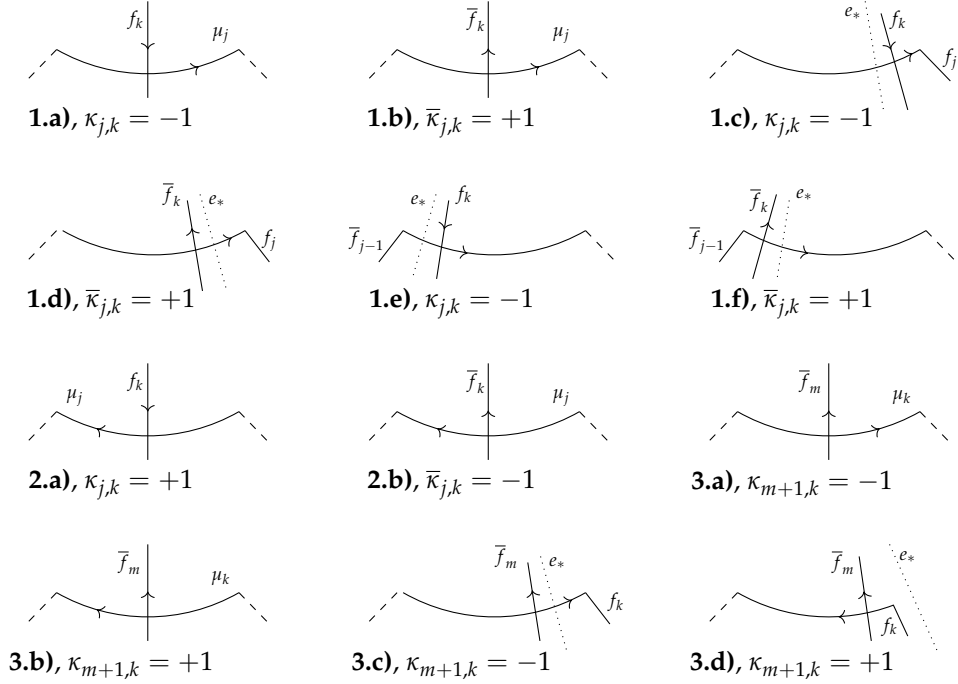
\end{proof}

\begin{remark} 
 The algorithm which we describe in the proof of Theorem \ref{theorem_algorithm} depends on certain choices, e.g. choosing the numbers $\epsilon_1,\ldots, \epsilon_m$ as well as the radii $r_2,\ldots, r_m$ to be increasing.
 The algorithm certainly does not depend on the specific values of the $\epsilon_i$ and of the radii.
    Moreover, we could have chosen to always \emph{go in mathematically positive direction} for the paths $\Xi_j$ or we could have made any other choice for the orientation of the circle segments.
    Any choice would have yielded an algorithm similar to the one which we have described above.
    The  sets $\mathcal{Q}$ and $\overline{\mathcal{Q}}$ do depend on these choices.
    While we do not claim that our algorithm is optimal in any sense, it seems quite intuitive to choose the paths $\Xi_j$ the way we did.
    If we had e.g. chosen to always go in the mathematical positive direction one would create many unnecessary self-intersections in certain situations.
    Moreover, we believe that the choices we made result in an algorithm which can relatively easily be made into an actual computer algorithm.
    It is furthemore clear that the string topology coproduct of $\Delta[\gamma]$ does not depend on any of these choices and hence any such algorithm gives the coproduct $\vee \Delta[\gamma]$.
\end{remark}

\begin{remark}
Recall that with integer coefficients a class of the form $\Delta[\gamma]\in {H}_1(LM,M)$ is not necessarily a generator.
In Section \ref{sec homology} we had introduced the number $l([\gamma])\in\mathbb{N}$ which satisfies that
$    \Delta[\gamma] = l([\gamma])\cdot [\widetilde{\gamma}]    $
for a generator $[\widetilde{\gamma}]\in H_1(LM,M)$.
  We note that Theorem \ref{theorem_algorithm} gives a computation of the string topology coproduct also for the generators $[\widetilde{\gamma}]$.
   By our construction of $M$ there is a cyclic word $v = c_{i_1}^{\pm}\ldots c_{i_m}^{\pm}$ whose associated loop is freely homotopic to $\gamma$.
    Theorem \ref{theorem_algorithm} thus gives us a computation of $\vee \Delta[\gamma]$.
    By linearity we must have that $\vee \Delta [\gamma]$ is divisble by $l([\gamma])$ and thus
    $$      \vee [\widetilde{\gamma}] =  \frac{1}{l([\gamma])} \Big( \vee \Delta [\gamma]    \Big) .     $$
    
    Recall further the computation of the homology of $LM$ from Section \ref{sec homology}.
    For degree reasons the coproduct can only be non-trivial in degree $1$.
    Thus Theorem \ref{theorem_algorithm} implicitly gives a complete computation of the coproduct for surfaces of genus $g\geq 2$.  
\end{remark}

\begin{example}\label{example_coproduct_algorithm}
    We illustrate the construction of the representative $\gamma$ and the computation of the coproduct for the case $g = 3$ and the word $v = c_4 c_6 c_3 c_1^{-1} c_5^{-1} c_4$.
    See Figure \ref{fig_example} for the respective steps of constructing $\gamma$.
     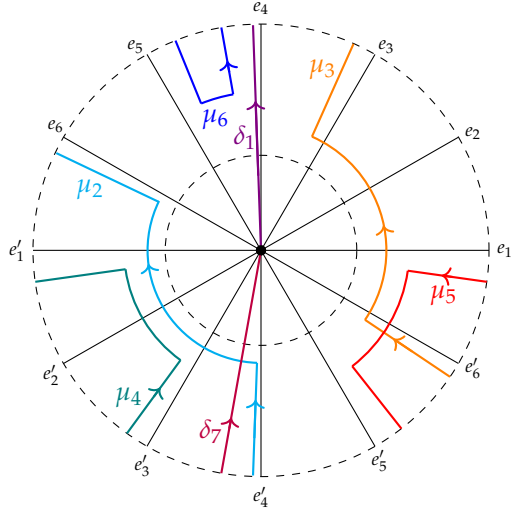
\begin{figure}
	\centering
	\begin{tikzpicture}[scale=3.]
             \filldraw (0.,0) circle (0.02cm);
            \draw[dashed] (0,0) circle (1.cm);
            \draw[dashed] (0.,0) circle (0.42cm);

            \draw[-] (1.,0) -- (-1,0) ;
            \draw[-] (-0.866,-0.5) -- (0.866,0.5);
             \draw[-] (-0.5,-0.866) -- (0.5,0.866);
            \draw[-] (0.,-1) -- (0,1);
            \draw[-] (-0.5,0.866) -- (0.5,-0.866);
            \draw[-] (-0.866,0.5) -- (0.866,-0.5);

            \draw (1.07,0) node[scale=0.7]{$e_1$};
            \draw (0.93,0.52) node[scale = 0.7]{$e_2$};
            \draw (0.55,0.9) node[scale = 0.7]{$e_3$};
            \draw (0.,1.07) node[scale = 0.7]{$e_4$};
            \draw (-0.55,0.9) node[scale = 0.7]{$e_5$};
            \draw (-0.9,0.55) node[scale = 0.7]{$e_6$};
            \draw (-1.07,0) node[scale=0.7]{$e_1'$};
            \draw (-0.92,-0.55) node[scale = 0.7]{$e_2'$};
            \draw (-0.52,-0.95) node[scale = 0.7]{$e_3'$};
            \draw (-0.,-1.09) node[scale = 0.7]{$e_4'$};
            \draw (0.52,-0.93) node[scale = 0.7]{$e_5'$};
            \draw (0.93,-0.52) node[scale = 0.7]{$e_6'$};

            \draw[-, thick ,violet] (0,0) -- (-0.035,0.995);
             \draw [->, thick ,violet] (-0.012,0.33) -- (-0.024,0.66);
     \draw[violet] (-0.07,0.5) node[scale=0.99]{$\delta_1$} ;

            \draw[-, thick ,cyan] (-0.035,-0.995) -- (-0.0175,-0.4975)    ; 
             \draw[-, thick, cyan] (154:0.4975) arc (154:268:0.4975) ;
            \draw[-,thick, cyan] (-0.898,0.43) -- (-0.45,0.215);
            \draw [->, thick, cyan] (-0.029,-0.83) -- (-0.0236,-0.66);
             \draw [->, thick, cyan] (-0.492,-0.07) -- (-0.494,-0.06);
        \draw[cyan] (-0.75,0.27) node[scale=0.99]{$\mu_2$} ;

            \draw[-, thick, orange] (0.829,-0.559) -- (0.456,-0.307);
            \draw[-, thick, orange] (-34:0.55) arc (-34:66:0.55) ;
            \draw[-, thick, orange] (0.406,0.913) -- (0.2233, 0.502);
            \draw [->, thick, orange] (0.663,-0.4472) -- (0.58,-0.391);
            \draw [->, thick, orange] (0.539,0.1) -- (0.537,0.11);
            \draw[orange] (0.27,0.8) node[scale=0.99]{$\mu_3$} ;

            \draw[-, thick, teal] (-0.588 ,-0.809) -- (-0.3528, -0.485);
             \draw[-, thick, teal] (234:0.6) arc (234:188:0.6) ;
             \draw[-, thick, teal] (-0.99,-0.139) -- ( -0.594,-0.0834);
             \draw[->, thick, teal] (-0.4704 ,-0.6472) -- (-0.441, -0.60675);
             \draw[teal] (-0.58,-0.65) node[scale=0.99]{$\mu_4$} ;

            \draw[-, thick, red] (0.99,-0.139) -- (0.6435, -0.09);
            \draw[-, thick, red] (308:0.65) arc (308:352:0.65) ;
            \draw[-, thick, red] ( 0.616,-0.788) -- (0.4, -0.512);
            \draw[->, thick, red] (0.8415, -0.118) -- ( 0.792,-0.111);
            \draw[red] (0.8,-0.2) node[scale=0.99]{$\mu_5$} ;

            \draw[-, thick, blue] (-0.375, 0.927) -- (-0.2625, 0.649);
            \draw[-, thick, blue] (112:0.7) arc (112:100:0.7) ;
            \draw[-, thick, blue] (-0.174,0.985) -- (-0.122, 0.69);
            \draw[->, thick, blue]  ( -0.139,0.788) --  (-0.148, 0.837);
            \draw[blue] (-0.2,0.58) node[scale=0.99]{$\mu_6$} ;

            \draw[-, thick, purple] (-0.174,-0.985) -- (0,0);
            \draw[->, thick, purple] (-0.1392, -0.788) -- ( -0.13 ,-0.739);
            \draw[purple] (-0.22,-0.8) node[scale=0.99]{$\delta_7$} ;
        
        \filldraw (0.,0) circle (0.02cm);
        
	\end{tikzpicture}

        \caption{The representative $\gamma$ for the word $c_4 c_6 c_3 c_1^{-1} c_5^{-1} c_4$. The parts $\mu_3$ and $\mu_5$ intersect, see case 2.b) in the proof of Theorem \ref{theorem_algorithm} as well as $\mu_2$ and $\delta_7$ which is case 3.d).}
        
        \label{fig_example}
    \end{figure}
    Note that at steps $2,3$ and $4$ we do not get any intersections.
    At step $5$ we have $s_5 = -1$ and since $ (f_5, \overline{f}_2, \overline{f}_4)= (e_5',e_6',e_1)\in\mathcal{C}_{4g}$ we have $\overline{\mathcal{S}}_5 = \{(5,2)\}$ and the sign is $\overline{\kappa}_{5,2} = -1$.
    This is case 2.b) and this intersection can be seen in Figure \ref{fig_example} as the intersection of $\mu_2$ and $\mu_5$.
    At step $6$ there is no intersection.
    In the last step, we obtain one more intersection which is case 3.d), since for $k = 2$ we have $s_k = -1$ and $\overline{f}_{{1}} = e_4'= \overline{f}_6$.
    Consequently, we obtain
    \begin{eqnarray*}
          \vee\Delta[\gamma] =   - [c_3 c_1^{-1}]\times [c_5^{-1} c_4 c_4 c_6] + [c_5^{-1} c_4 c_4 c_6]\times [ c_3 c_1^{-1}]  \\  +\, [c_6 c_3 c_1^{-1} c_5^{-1} c_4]\times [ c_4]  - [c_4]\times [ c_6 c_3 c_1^{-1} c_5^{-1} c_4]   . 
    \end{eqnarray*}
    
\end{example}

\section{String cobracket and Turaev Cobracket}
\label{subsec Turaev cobracket}

In this section we compare the string cobracket on surfaces to the Turaev cobracket.
Our definition of the cobracket follows \cite{naefwillwacher} where the cobracket is defined on the equivariant homology of the free loop space, in contrast to \cite{GoreskyHingston} who define the cobracket on relative equivariant homology.
We only consider homology with rational coefficients in this section.

We begin by defining the Turaev cobracket.
Let again $L$ be the set of free homotopy classes on a surface of genus $g \geq 2$, and let $V$ be the $\mathbb{Q}$-vector space spanned by $L$.
Further, let $\widetilde{L}\subseteq L$ be the set of non-trivial free homotopy classes and $\widetilde{V}$ the corresponding $\mathbb{Q}$-vector space spanned by $\widetilde{L}$.
Canonically, the space $\widetilde{V}$ is a subspace of $V$ and there is a canonical projection map $p\colon V\to \widetilde{V}$ with kernel $\mathbb{Q}[\gamma_0]$, where $[\gamma_0]$ is the class of the contractible loop.
The Turaev cobracket is defined as follows, see Section 8.1 in \cite{Turaev}.

\begin{definition} 
    Let $\gamma$ be a loop such that each self-intersection is a transversal double point $q$, i.e. the two tangent vectors $v_1, v_2$ of $\gamma$ at a self-intersection point $q$ are linearly independent. We choose the order such that $(v_1, v_2)$ is positively oriented. 
    Denote by $\gamma_i^q$ the part of $\gamma$ starting at $q$ in the direction of $v_i$ until it reaches $q$ again. The Turaev cobracket $\vee_T: V \xrightarrow[]{} V \otimes V$ is defined on the generators as  
    \begin{equation}
        \vee_T([\gamma]) = \sum_{q} \big(p ([\gamma_1^q]) \otimes p([\gamma_2^q]) - p([\gamma_2^q] ) \otimes p([\gamma_1^q])\big)
    \end{equation}
    where the sum runs over all intersection points $q$.
\end{definition}

\noindent Turaev shows that $\vee_T$ defines a Lie cobracket on $V$, and together with the Goldman Lie bracket makes $V$ into a Lie bialgebra, see \cite{Turaev}, Theorem 8.3. 
We want to compare the Turaev cobracket to the string cobracket which we now introduce.

Let $M$ be a topological space and consider the free loop space $LM = \mathrm{Map}({S}^1,M)$.
The free loop space has a circle action ${S}^1\times LM\to LM$ which is not free, since the constant loops are always fixed and also powers of loops have non-trivial stabilizers.
Recall that there is a universal ${S}^1$-space $E{S}^1$ which is a contractible space and which has a free ${S}^1$-action.
We consider the Borel construction $$ \pi\colon LM \times E{S}^1 \to LM_{{S}^1} :=  (LM\times E{S}^1)/{S}^1 $$
where we quotient out the diagonal action which is free because of the freeness of the $S^1$-action on $ES^1$.
A free ${S}^1$-action induces a principal bundle, hence the map $\pi$ is a circle bundle.
Moreover, since $E{S}^1\simeq \{\mathrm{pt}\}$ we have $LM\times E{S}^1\simeq LM$ and the Gysin sequence of $\pi$ thus induces a long exact sequence of the form
$$      \ldots \xrightarrow[]{} {H}_{p+2}(LM ) \xrightarrow[]{\pi_*} {H}_{p+2}(LM_{{S}^1}) \xrightarrow[]{\chi\cap} {H}_p(LM_{{S}^1}) \xrightarrow[]{\sigma_*} {H}_{p+1} (LM) \xrightarrow[]{} \ldots      $$
where $\chi\in{H}^2(LM_{S^1})$ is the Euler class and the map $\sigma$ is induced by the transfer morphism of the circle bundle $\pi$.
By definition the homology of the homotopy quotient $LM_{{S}^1}$ is the ${S}^1$-equivariant homology of $LM$.
We therefore get maps
$$    \mu\colon  {H}_i^{{S}^1}( LM) \xrightarrow[]{\sigma_*} {H}_{i+1}(LM) \xrightarrow[]{} {H}_{i+1}(LM,M)     $$
as well as
$$    \eta \colon {H}_j(LM,M) \xrightarrow[]{} {H}_j(LM) \xrightarrow[]{\pi_*}  {H}_j^{{S}^1}(LM)      $$
where the map $H_j(LM,M)\to H_j(LM)$ is induced by the long exact sequence of the pair $(LM, M)$ which is split since the inclusion $M\hookrightarrow LM$ has the evaluation as a left-inverse.
The \emph{string cobracket} is defined as the composition
$$  \vee^{S^1}\colon   H_i^{S^1}(LM)    \xrightarrow[]{\mu}  H_{i+1}(LM,M) \xrightarrow[]{\vee} \big(H_{*}(LM,M)^{\otimes 2} \big)_{i+2-n}\xrightarrow[]{\eta\otimes \eta} \big(H_{*}^{S^1}(LM)^{\otimes 2}\big)_{i+2-n}. $$
Note that for $M$ a closed oriented surface of genus $g\geq 2$ the cobracket is trivial for degree reasons for $i\neq 0$.
Therefore we will only consider the equivariant homology of degree $0$ from now on.

The connected components of $LM_{{S}^1}$ are given by the images of the connected components of $LM\times E{S}^1$ under the quotient map $\pi\colon LM\times E{S^1}\to LM_{{S}^1}$.
Hence, the components of $LM_{{S}^1}$ are in one-to-one corresponendence with the conjugacy classes of $\pi_1(M)$.
Choose a basepoint $u_0\in E{S}^1$.
There is a map $\psi \colon V\to H_0^{S^1}(LM)$ induced by mapping a free homotopy class $[\gamma]\in \widetilde{L}$ to the class $[\pi(\gamma,u_0)]\in H_0^{S^1}(LM)$.
This map is clearly an isomorphism.
We can thus consider the Turaev cobracket as a map $\vee_T \colon H_0^{S^1}(LM) \to H_0^{S^1}(LM)^{\otimes 2}$.
We first need an auxilliary lemma.

\begin{lemma}\label{lemma_transfer_map}
    Let $[h]$ be a conjugacy class in $\pi_1(M)$ and let $\gamma\in LM$ be a representative.
    Then 
    $$    \mu([\pi(\gamma,u_0)] ) =  \Delta ([\gamma])          $$
    where $\Delta \colon H_0(LM)\to H_1(LM)$ is the BV operator as before.
\end{lemma}
\begin{proof}
    The map $\sigma_*\colon {H}_0(LM_{{S}^1})\to {H}_1(LM)$ is given as the composition
    $$   {H}_0(LM_{{S}^1}) \xrightarrow[]{\sigma'_*} {H}_1(LM\times E{S}^1)     \xrightarrow[\cong]{(\mathrm{pr}_1)_*}  {H}_1(LM)         $$
    where $\sigma'$ is the transfer morphism of the Gysin sequence for the circle bundle $LM\times ES^1\to LM_{S^1}$.
    Consider the general situation of a $(k-1)$-sphere bundle $\pi\colon SE \to B$ over a path-connected space $B$. 
    The transfer map $\sigma$ is obtained by the upper part of the commuting diagram
    $$
    \begin{tikzcd}
        {H}_0(B) \arrow[]{rd}{\sigma_*} & 
        \\
        {H}_{k}(DE,SE) \arrow[]{u}{\pi_*(\tau\cap (\cdot))}
        \arrow[hook]{r}{\partial}
        & 
        {H}_{k-1}(SE) 
        \\
        {H}_k(D^k,{S}^{k-1}) \arrow[]{u}{s_*} \arrow[]{r}{\partial,\cong}
        &
        {H}_{k-1}({S}^{k-1}) \arrow[]{u}{s_*}
    \end{tikzcd}
    $$
    where $s\colon (D^k,{S}^{k-1}) \hookrightarrow (DE,SE)$ is the inclusion of the fiber pair.
    Let $p_0\in B$ be a point and consider the corresponding generator $[p_0]\in H_0(B)$.
    We claim that $\sigma_*[p_0] = s_* [S^{k-1}]$ where $[S^{k-1}]\in H_{k-1}(S^{k-1})$ is the positively oriented generator with respect to the standard orientation of the sphere.
    Let $[D^k]\in H_k(D^k,S^{k-1})$ be the relative homology class which is mapped to $[S^{k-1}]$ under the connecting morphism $\partial\colon H_k(D^k,S^{k-1})\to H_{k-1}(S^{k-1})$.
    By the construction of the Thom class $\tau\in \mathrm{H}^k(DE,SE)$ we have that $\langle s^*\tau, [D^k]\rangle = 1$, where $\langle\cdot,\cdot\rangle$ is the Kronecker pairing.
    Note that the by the properties of the Kronecker pairing it holds that $\langle s^*\tau ,[D^k]\rangle = \langle \tau, s_*[D^k]\rangle$.
    Moreover, if $[q_0]\in H_0(DE)$ is the generator in degree $0$ induced by a point $q_0\in DE$, the Kronecker pairing satisfies
    $$ \tau\cap (s_*[D^k]) =  \langle \tau, s_*[D^k]\rangle \cdot [q_0]  = [q_0] .      $$
    Since the class $[q_0]\in H_0(E)$ is mapped to the class $[p_0]\in H_0(B)$ under the map $\pi_*$ we obtain that $\pi_*(\tau\cap (s_*[D^k])) = [p_0]$.
    It follows from the commutative diagram that $\sigma_*[p_0] = s_* [S^{k-1}]$.
    
    Specializing to our present situation we have that the inclusion of the fiber $s'\colon {S}^1\to LM\times E{S}^1$ is the orbit map induced by the ${S}^1$-action.
    Hence under the projection on the first factor $\mathrm{pr}_1\colon LM\times E{S}^1\to LM$ we see that $\mathrm{pr}_1 \circ s'\colon {S}^1\to LM$ is the orbit of the usual ${S}^1$-action on $LM$ and thus $  \mu([\pi(\gamma,u_0)] ) =  \Delta ([\gamma]) $.
\end{proof}

\begin{theorem}\label{theorem_turaev}
    Let $M$ be a closed oriented surface of genus $g\geq 2$.
    Under the canonical identification $V \cong {H}_0^{{S}^1}(LM)$ the Turaev cobracket and the negative of the string cobracket agree, i.e. the diagram
    $$
\begin{tikzcd}
    {H}_1(LM,M) \arrow[]{r}{-\vee} & {H}_0(LM,M)^{\otimes 2} \arrow[]{d}{\eta\otimes \eta}
    \\
    {H}_0^{{S}^1}(LM) \arrow[swap]{u}{\mu} \arrow[]{r}{\vee_T} & {H}_0^{{S}^1}(LM)^{\otimes 2} 
\end{tikzcd}
$$
commutes.
\end{theorem}
\begin{proof}
    Let $[\gamma]\in V$ be a free homotopy class of loops in $M$.
    By Lemma \ref{lemma_transfer_map} we need to compare $\vee_T [\gamma]$ and $(\eta\otimes \eta )  \big( \vee \Delta[\gamma] \big)  $.
    We choose a representative $\gamma$ with only transverse self-intersections, e.g. by applying Theorem \ref{theorem_algorithm} if desired.
    Lemma \ref{lemma_coproduct_signs} yields the coproduct.
    One sees that the signs that we obtain from Lemma \ref{lemma_coproduct_signs} are precisely minus the signs that one gets in the definition of the Turaev cobracket.
\end{proof}

\begin{remark}
    We note that the following change in the definition of the coproduct would yield the agreement of the Turaev cobracket with the string cobracket instead of with the negative of the string cobracket as in the above theorem.
    As we have seen in Section \ref{sec string operations} when computing the coproduct, we take the homology cross product of a class $X\in H_*(LM,M)$ with the fundamental class $[I]\in H_1(I,\partial I)$ of the unit interval, i.e. mapping $X\mapsto X\times [I]$.
    If we instead swap the factors, i.e. we take the cross product $X\mapsto [I]\times X$ then for classes $X$ of odd degree the coproduct of $X$ gets an additional minus sign.
    In particular the proof of Theorem \ref{theorem_turaev} shows that we would then get an agreement of the Turaev cobracket with the string cobracket.

    We further note that the choice of the identification $\Phi\colon N(\Delta M)\to TM$ as in Section \ref{sec string operations} of the normal bundle $N(\Delta M)$ with the tangent bundle $TM$ does not influence the signs of the coproduct and cobracket in even dimensions.
\end{remark}

\section{The coproduct for hyperbolic manifolds and $3$-manifolds}

After having studied the string topology coproduct and cobracket on surfaces, it seems natural to consider the following two situations.
After considering $2$-dimensional manifolds the next step is to ask if certain properties carry over to $3$-dimensional manifolds.
Moreover, since all orientable surfaces of genus $g\geq 2$ are hyperbolic manifolds, one can ask about the coproduct on arbitrary hyperbolic manifolds.
In this section we give some remarks on how the string topology coproduct behaves for hyperbolic manifolds in higher dimensions and what can be said about the string topology coproduct on $3$-manifolds.

\begin{theorem}\label{theorem_triviality_cyclic_centralizers}
    Let $M$ be an aspherical manifold of dimension $n\geq 3$ with the property that for each $h\in\pi_1(M)$, $h\neq e$ the centralizer $C_h$ has homological dimension less than $ n -1$.
    Then the string topology coproduct on $M$ is trivial.
\end{theorem}
\begin{proof}
    Let $[h]$ be a conjugacy class in $\pi_1(M)$ and let $\gamma\in LM$ be a representative of $[h]$.
    As in Section \ref{sec homology}, one sees directly that the component of $LM$ containing $\gamma$ is an Eilenberg-MacLane space $K(C_h,1)$ where $C_h$ is the centralizer of $h$ in $\pi_1(M)$.
    By assumption we have ${H}_i(LM_{[h]}) \cong \{0\}$ for $i\geq n$.
    Consequently, every non-trivial class $X\in{H}_{*}(L M,M)$ satisfies $\mathrm{deg}(X) < n-1$.
    The coproduct is thus trivial for degree reasons.
\end{proof}
\begin{cor}
    Let $M$ be a closed hyperbolic manifold of dimension $n\geq 3$.
    Then the string topology coproduct of $M$ vanishes.
\end{cor}
\begin{proof}
    As pointed out in the discussion before Lemma \ref{lemma_centralizers_hyperbolic_grps}, the centralizers of non-trivial elements in the fundamental group of a hyperbolic manifold are infinite cyclic and hence the above theorem applies.
\end{proof}
\begin{cor}\label{cor_coproduct_trivial_alg_hyp_3_mfld}
    Let $M$ be a closed aspherical $3$-manifold such that $\pi_1(M)$ is algebraically hyperbolic, i.e. the fundamental group has no rank $2$ abelian subgroup.
    Then the coproduct of $M$ and of all its finite covering spaces vanishes.
\end{cor}
\begin{proof}
    By \cite[Lemmas 11.2 and 11.3]{abbaspour2005string} the assumptions of Theorem \ref{theorem_triviality_cyclic_centralizers} are met. 
\end{proof}
Abbaspour shows in \cite{abbaspour2005string} that there is a dichotomoy on $3$-manifolds concerning the Chas-Sullivan product.
Note that there is a canonical splitting ${H}_{*}(L M)\cong {H}_{*}(L M,M)\oplus {H}_{*}(M)$.
Abbaspour shows that if $M$ is a $3$-manifold which is not algebraically hyperbolic 
then the restriction of the loop product to ${H}_{*}(L M,M)$ is non-trivial for $M$ or a double cover of $M$.
On the other hand if the manifold is algebraically hyperbolic, then the restriction of the Chas-Sullivan product to ${H}_{*}(L M,M)$ is trivial.
As we have seen in Corollary \ref{cor_coproduct_trivial_alg_hyp_3_mfld}, the coproduct is trivial for algebraically hyperbolic $3$-manifolds.
It is therefore natural to ask whether the string topology coproduct also detects the dichotomy between algebraically hyperbolic $3$-manifold and those which are not algebraically hyperbolic.
This would mean that the string coproduct should be non-trivial if and only if the manifold is not algebraically hyperbolic.
Note that the string topology coproduct is known to be non-trivial for the $3$-sphere, see \cite{HingstonWahl}.
However, as shown in \cite{Stegemeyer} the string topology coproduct is trivial for the $3$-dimensional torus.
Hence, the dichotomy that holds for the Chas-Sullivan product, see \cite{abbaspour2005string}, does not hold analogously for the string topology coproduct.


\end{document}